\documentclass[final,leqno,onefignum,onetabnum]{article}

\usepackage{amsmath}
\usepackage{amsfonts}
\usepackage{epsfig}
\usepackage{multicol}
\usepackage{graphicx, graphics}
\usepackage{color}

\newtheorem{thm}{Theorem}[section]
\newtheorem{proposition}[thm]{Proposition}
\newtheorem{lemma}[thm]{Lemma}

\newcommand{\eop}{{\hspace*{\fill}$\square$}}

\title{An NPZ Model with State-Dependent Delay due to Size-Structure in Juvenile Zooplankton}

\author{Matt Kloosterman\footnote{mklooste@uwaterloo.ca} \and
        Sue Ann Campbell \and
        Francis J. Poulin \\ Department of Applied Mathematics, \\University of Waterloo, Waterloo, ON, Canada N2L 3G1 
}

\begin{document}
\maketitle

\begin{abstract}
The study of planktonic ecosystems is important as they make up the bottom trophic levels of aquatic food webs. We study a closed Nutrient-Phytoplankton-Zooplankton (NPZ) model that includes size structure in the juvenile zooplankton. The closed nature of the system allows the formulation of a conservation law of biomass that governs the system. The model consists of a system of nonlinear ordinary differential equation coupled to a partial differential equation. We are able to transform this system into a system of delay differential equations where the delay is of threshold type and is state-dependent. The system of delay differential equations can be further transformed into one with fixed delay. Using the different forms of the model we perform a qualitative analysis of the solutions, which includes studying existence and uniqueness, positivity and boundedness, local and global stability, and conditions for extinction. Key parameters that are explored are the total biomass in the system and the maturity level at which the juvenile zooplankton reach maturity. Numerical simulations are also performed to verify our analytical results.
\end{abstract}

\section{Introduction}

Nutrient-Phytoplankton-Zooplankton (NPZ) models are used to describe the bottom two trophic levels of an aquatic ecosystem. As is the case with many ecological models, they range from very simple to very complex. Simple models, such as the Lotka-Volterra system \cite{murray}, are beneficial in that one can obtain analytical results more easily, but often suffer from lack of realism. On the other hand, complex models may theoretically represent a more accurate description of reality, but may be difficult or impossible to understand in any general way, and may be useless without precise and accurate parameter values. A lot of structure can be given to an NPZ model, such as size-dependence and spatial dependence, and this structure leads to rich modelling possibilities such as size-dependent parameters and coupling with external factors like fluid dynamics or higher predation. To strike a balance, it may be useful to start with a simple model that focuses on one factor affecting the ecosystem, and use it to obtain analytical results and study general trends. For an introduction in the construction of NPZ models, see \cite{franks2002}.

We will focus on the role of maturity in the juvenile zooplankton population within an ecosystem while keeping all other factors as simple as possible. In this vain, a simple NPZ model is coupled with a standard linear first-order PDE that describes the spectrum of the juvenile zooplankton population as a function of time and maturity. It is a known result that this type of PDE is closely related to delay equations. Examples of population models where time delay is a consequence of age structure can be found in \cite{gurney1983}, \cite{blythe1984}, \cite{sulsky1989}, \cite{bocharov2000}. We will be considering size structure where the rate of growth of the juvenile zooplankton is permitted to depend on the concentration of the phytoplankton. With this assumption, the delay in the related delay equation is not explicitly defined, but rather defined implicitly through a threshold-type condition. Consequently, systems with this type of delay are known as threshold delay differential equations. They have been studied in the case of a single population in \cite{smith1993} and for an insect species in \cite{nisbet1983}, for example. These models operate under the assumption that maturation occurs when an immature individual accumulates enough of some quantity, such as size or weight. Other applications include red blood cell production \cite{metz} and immunology \cite{waltman1977}.

Our model is formulated in such a way that biomass is conserved. In other words, we are considering a closed ecosystem with no mass being added to or subtracted from the system. NPZ models with this property have been studied, for example, by \cite{gentleman2008} and \cite{wroblewski1988}. A typical property of these systems is that the amount of biomass, which is determined by initial conditions, plays a crucial role in determining the types of dynamics that can and do occur. Most notably, an insufficient amount of biomass leads to the extinction of plankton; either the zooplankton only, or both the phytoplankton and zooplankton.

\section{Structured Model}

We consider a Nutrient-Phytoplankton-Zooplankton (NPZ) model in which the zooplankton population is split into a mature class and a juvenile class. We will consider the juvenile class as a function of time and maturity, denoted $\rho(t,s)$. Maturity is considered to be an abstract quantity, of which the juvenile zooplankton much accumulate enough in order to enter adulthood. The total concentration of immature zooplankton with maturity levels between $s_1$ and $s_2$ at time $t$ is then $\int_{s_1}^{s_2} \rho(t,s) \,ds$. We restrict $s$ to the interval $[0,m]$, where $m$ is the required level of maturity for adulthood. The nutrient, phytoplankton, and mature zooplankton variables will therefore only depend on time.

Since the zooplankton feed on phytoplankton, we will assume that the growth rate of the juvenile zooplankton depends on the quantity of phytoplankton present in the system. We will denote this dependence as $R(P(t))$.

The following equations model the ecosystem and we will refer to it as the PDE1 model.
\begin{subequations}
\label{model_PDE}
\begin{align}
&\frac{dN(t)}{dt} = -\mu P(t)f(N(t)) +\lambda P(t) + \delta Z(t) + (1-\gamma)gZ(t)h(P(t)) \nonumber\\
&\hspace{1cm}+ \int_0^m \delta_0 \rho(t,s) \,ds, \label{model_PDE_N}\\
&\frac{dP(t)}{dt} = \mu P(t)f(N(t)) - \lambda P(t) - g Z(t) h(P(t)), \label{model_PDE_P}\\
&\frac{dZ(t)}{dt} = R(P(t)) \rho(t,m) - \delta Z(t), \label{model_PDE_Z}\\
&\frac{\partial \rho}{\partial t}(t,s)+R(P(t))\frac{\partial \rho}{\partial s}(t,s) = -\delta_0 \rho(t,s), \label{model_PDE_rho}\\
&R(P(t))\rho(t,0)=\gamma g Z(t) h(P(t)). \label{model_PDE_BC}
\end{align}
\end{subequations}
Appropriate initial conditions for these equations are
\begin{align}
\label{model_PDE_IC}
    N(0) = N_0, \hspace{0.5cm} P(0)=P_0, \hspace{0.5cm}Z(0) = Z_0,\hspace{0.5cm} \rho(0,s) = \rho_0(s),
\end{align}
where $N_0, P_0,$ and $Z_0$ are nonnegative real numbers and $\rho_0$ is a nonnegative function on the interval $[0,m]$.

A few basic ecological processes govern the system. Phytoplankton, $P$, uptake nutrient, $N$, at a rate that is proportional to the quantity of phytoplankton and a function of the dissolved nutrient, $f(N)$.  Zooplankton, $Z$, graze on phytoplankton at a rate that is proportional to the quantity of zooplankton and a function of  phytoplankton, $h(P)$. This is marked by a grazing efficiency factor, $\gamma \in (0,1]$. Phytoplankton mortality is proportional to the amount of phytoplankton present at the time. Zooplankton mortality is similarly proportional to the amount of zooplankton in the system at that time. Nutrient recycling transforms dead biomass and faecal matter back to the dissolved nutrient variable. Equation \eqref{model_PDE_rho} is a standard transport equation \cite{metz} with a decay term due to natural mortality. The boundary condition, given in equation \eqref{model_PDE_BC}, says that the birth rate of the immature zooplankton is equal to some fraction of the biomass obtained by the zooplankton through grazing.

Some assumptions are made on the functional responses. For the phytoplankton nutrient uptake response, we assume $f\in \mathcal{C}^1$ and
\begin{align}
\label{f_props}
    f(0)=0, \hspace{0.5cm} f'(N) >0, \hspace{0.5cm}  \lim_{N\rightarrow \infty} f(N)=1.
\end{align}
Similarly, for the zooplankton functional response for grazing on phytoplankton, we assume $h\in \mathcal{C}^1$ and
\begin{align}
\label{h_props1}
    h(0)=0, \hspace{0.5cm} h'(P) > 0, \hspace{0.5cm}  \lim_{P\rightarrow \infty} h(P)=1.
\end{align}
These assumptions on $h$ encompass both Type II (concave down) and Type III (sigmoidal shaped) responses, as described in \cite{holling}.

We will also assume throughout that all parameter values are positive (although we will allow $\delta_0$ to be zero) and that
\begin{align*}
    \mu>\lambda, \hspace{0.5cm} \gamma g > \delta,
\end{align*}
so that $f^{-1}(\lambda/\mu)$ and $h^{-1}(\delta/(\gamma g))$ exist.

It will be assumed that $R\in \mathcal{C}^1$ and that it satisfies the following properties.
\begin{align}
\label{Rprops}
    R(P)\geq 0,\hspace{0.2cm} R'(P)\geq 0, \hspace{0.2cm}R'(0)>0 \mbox{ if } R(0)=0, \hspace{0.2cm} \lim_{P\rightarrow \infty} R(P)=R_{\infty}< \infty.
\end{align}
These assumptions can be justified by laboratory experiments in \cite{mccauley}, where they measured the development rate of the zooplankter \emph{Daphnia} as a function of $F$, the amount of available food. They found that the development rate was proportional to $F/(F+F_{half})$.

The assumptions on $h$ in equation \eqref{h_props1} and the third assumption on $R$  in equation \eqref{Rprops} imply that
\begin{align}
\label{hR_limit}
    \lim_{P\rightarrow 0^+} \frac{h(P)}{R(P)} < \infty.
\end{align}
This puts a bound on $\rho(t,0)$ as the phytoplankton population approaches zero, as seen in equation \eqref{model_PDE_BC}.

As previously stated, there is no biomass lost from system \eqref{model_PDE}. We can confirm this by noting that PDE1 model \eqref{model_PDE} satisfies the following conservation law:
\begin{align}
\label{conlaw_PDE}
    N_T = N(t) +P(t) +Z(t) + \int_0^m \rho(t,s)\,ds,
\end{align}
where $N_T$ is the total biomass in the system. Given this conservation law, the total amount of juvenile zooplankton can be determined from the total biomass and the quantity of the dissolved nutrient, phytoplankton, and zooplankton. Making the appropriate substitution of equation \eqref{conlaw_PDE} into \eqref{model_PDE} yields the following system:
\begin{subequations}
\label{model_PDE_2}
\begin{align}
&\frac{dN(t)}{dt} = -\mu P(t)f(N(t)) +\lambda P(t) + \delta Z(t) + (1-\gamma)gZ(t)h(P(t)) \nonumber\\
&\hspace{2cm}+ \delta_0(N_T-N(t)-P(t)-Z(t)), \label{model_PDE_N_2}\\
&\frac{dP(t)}{dt} = \mu P(t)f(N(t)) - \lambda P(t) - g Z(t) h(P(t)), \label{model_PDE_P_2}\\
&\frac{dZ(t)}{dt} = R(P(t)) \rho(t,m) - \delta Z(t), \label{model_PDE_Z_2}\\
&\frac{\partial \rho}{\partial t}(t,s)+R(P(t))\frac{\partial \rho}{\partial s}(t,s) = -\delta_0 \rho(t,s), \label{model_PDE_rho_2}\\
&R(P(t))\rho(t,0)=\gamma g Z(t) h(P(t)). \label{model_PDE_BC_2}
\end{align}
\end{subequations}
We will refer to this system as the PDE2 model. While the substitution is straightforward, it is important to note that we have now fixed $N_T$. Thus, in order for a solution of system \eqref{model_PDE_2} to be the same as a solution to system \eqref{model_PDE}, we must choose initial conditions in \eqref{model_PDE_IC} that satisfy
\begin{align}
\label{conlaw_PDE_IC}
    N_T = N_0 +P_0 +Z_0 + \int_0^m \rho_0(s)\,ds.
\end{align}
In other words, the total biomass, $N_T$, is set by the initial conditions for the PDE1 model \eqref{model_PDE}, so if we treat it as a parameter in the PDE2 model \eqref{model_PDE_2}, then we have to restrict the possible initial conditions by the relation \eqref{conlaw_PDE_IC}. If we fix $N_T$ in system \eqref{model_PDE_2} and choose initial conditions that do not satisfy \eqref{conlaw_PDE_IC}, then the resulting solution is not a solution to PDE1. Indeed, if we define
\begin{align}
\label{Delta}
   \Delta(t) =-N_T + N(t) +P(t) +Z(t) + \int_0^m \rho(t,s)\,ds,
\end{align}
it can be verified that a solution to \eqref{model_PDE_2} satisfies
\begin{align}
\label{Delta_DE}
   &\frac{d}{dt} \Delta(t) = -\delta_0 \Delta(t).
\end{align}
This implies that if equation \eqref{conlaw_PDE_IC} is not satisfied ($\Delta(0)\neq 0$), then equation \eqref{conlaw_PDE} is never satisfied ($\Delta(t)\neq 0$). However, a solution for PDE1 always satisfies equation \eqref{conlaw_PDE} ($\Delta(t)\equiv0$). We can conclude, however, that a solution to system PDE2 will asymptotically approach a solution to PDE1, assuming that both solutions exist for all time. Due to the simpler nature of the PDE2 model, we will study it instead of the original system, and assume that initial conditions satisfy equation \eqref{conlaw_PDE_IC}.

\cite{smith1993} studies a model in a similar form to the PDE2 model \eqref{model_PDE_2}, although the system is for a single species where the mature population is coupled with its own immature class. The growth rate of the immature class depended on the size of the mature population. \cite{smith1993} reduced the equations so that they became a scalar threshold delay equation. This type of equation has a state-dependent delay due to the fact that a new-born individual must first reach maturity before entering into the adult population. This delay is defined implicitly via an integral which represents the condition for reaching maturity. \cite{smith1993} then uses a clever change of variables to transform the system to a functional differential equation where the delay no longer depends on the state. Due to the similar form of \eqref{model_PDE} to Smith's model, we will use the same methods to reduce our model to a threshold delay equation and then to a delay differential equation with a bounded and state-independent delay.

\section{Reduction to a Threshold Delay Equation}

As in \cite{smith1993}, we can solve \eqref{model_PDE_rho_2} and \eqref{model_PDE_BC_2} for $\rho$ as a function of $P$ and $Z$ using the method of characteristics, and then substitute this solution into \eqref{model_PDE_Z_2}. It can be verified that the solution for $\rho$ is
\begin{equation}
\label{rho_solution}
\rho(t,s) = \left\{ \begin{array}{ll}
    e^{-\delta_0 t}\rho_0\left(s-\int_{0}^t R(P(u))\,du\right) & \mbox{ for } (t,s) \in S_1 \\
    e^{-\delta_0 \tau(s,P_t)}\frac{\gamma g Z(t-\tau(s,P_t))h(P(t-\tau(s,P_t)))}{R(P(t-\tau(s,P_t)))} & \mbox{ for } (t,s) \in S_2
    \end{array} \right.,
\end{equation}
where $P_t(u) = P(t+u)$, and $\tau(s,P_t)$  is defined implicitly by
\begin{align}
\label{tau_equation}
    \int_{-\tau(s,P_t)}^0 R(P(t+u))\,du=s.
\end{align}
Biologically, it is the time it takes for the immature zooplankton to grow to size $s$ as a function of the history of phytoplankton. For example, if the phytoplankton population is at an equilibrium value $P(t)=P^*$ for all $t$, then $\tau(s,P_t)=\frac{s}{R(P^*)}$.

Substituting $\rho(t,m)$, defined by equation \eqref{rho_solution}, into \eqref{model_PDE_Z_2} gives a nonautonomous threshold delay equation. Alternatively, we can change our initial conditions from \eqref{model_PDE_IC} to
\begin{align}
\label{model_TDE_IC}
    N(t_0+t) = \phi_1(t), \hspace{0.2cm}  P(t_0+t)=\phi_2(t),\hspace{0.2cm}  Z(t_0+t)=\phi_3(t),
\end{align}
for $ t\in [-t_0,0]$, with $\tau(m,\phi_2) = t_0$. The initial conditions are shifted in a non-standard way so they correspond with the original PDE1 model \eqref{model_PDE}. The time $t_0$ is when all of the immature zooplankton initially present at $t=0$ have reached maturity. With these initial conditions, we can restrict $\rho(t,s)$ to $S_2$. Then we have the following autonomous system
\begin{subequations}
\label{model_TDE_2}
\begin{align}
\frac{dN(t)}{dt} =& -\mu P(t)f(N(t)) +\lambda P(t) + \delta Z(t) + (1-\gamma)gZ(t)h(P(t)) \nonumber\\
    &+ \delta_0(N_T-N(t)-P(t)-Z(t)), \label{model_TDE2_N}\\
\frac{dP(t)}{dt} =& \mu P(t)f(N(t)) - \lambda P(t) - g Z(t) h(P(t)), \label{model_TDE2_P}\\
\frac{dZ(t)}{dt} =& R(P(t)) e^{-\delta_0 \tau(m,P_t)}\frac{\gamma g Z(t-\tau(m,P_t))h(P(t-\tau(m,P_t)))}{R(P(t-\tau(m,P_t)))} - \delta Z(t), \label{model_TDE2_Z}
\end{align}
\end{subequations}
for $t\geq t_0$, where the function $\tau$ is given by equation \eqref{tau_equation}.  We will refer to this system as the TDE model, as systems in this form are known as threshold delay equations.

Substituting $\rho(t,s)$ for $(t,s)\in S_2$ into the definition of $\Delta$ in equation \eqref{Delta}, we find that
\begin{align}
\label{Delta_TDE}
    \Delta(t) =&-N_T + N(t) +P(t) +Z(t) \nonumber\\ &+ \int_0^m e^{-\delta_0 \tau(s,P_t)}\frac{\gamma g Z(t-\tau(s,P_t))h(P(t-\tau(s,P_t)))}{R(P(t-\tau(s,P_t)))}\,ds.
\end{align}
It can be verified that solutions to system \eqref{model_TDE_2} still satisfy equation \eqref{Delta_DE} with $\Delta$ defined by equation \eqref{Delta_TDE}. In particular, if the initial conditions \eqref{model_TDE_IC} satisfy
\begin{align}
\label{conlaw_IC}
    N_T =& \phi_1(0) + \phi_2(0) + \phi_3(0) \nonumber \\ &+ \int_0^m e^{-\delta_0 \tau(s,\phi_2)}\frac{\gamma g \phi_3(-\tau(s,\phi_2))h(\phi_2(-\tau(s,\phi_2)))}{R(\phi_2(-\tau(s,\phi_2)))} \,ds,
\end{align}
then we have
\begin{align}
\label{conlaw_TDE}
    N_T =& N(t) + P(t) + Z(t)\nonumber \\  &+ \int_0^m e^{-\delta_0 \tau(s,P_t)}\frac{\gamma g Z(t-\tau(s,P_t))h(P(t-\tau(s,P_t)))}{R(P(t-\tau(s,P_t)))} \,ds,
\end{align}
for $t\geq t_0$. Since initial conditions satisfying equation \eqref{conlaw_IC} are the only ones that relate to the original system \eqref{model_PDE}, we will restrict our attention to this particular choice. For later use we define the set
\begin{align*}
    D_{N_T} = &\{ (\phi_1,\phi_2,\phi_3)^T \in \mathcal{C}([-t_0,0],\mathbb{R}^3) : \phi_i(t) >0, i=1,2,3, t\in[-t_0,0], \\ & \tau(m,\phi_2) = t_0, \mbox{ and \eqref{conlaw_IC} is satisfied}  \}.
\end{align*}
This is the set of initial conditions that we will consider. That is, we will not study cases where $N,P,$ or $Z$ have initial values equal to zero as we are interested in studying the dynamics when all three trophic levels are initially present.

We note that system \eqref{model_TDE_2} is a differential equation with locally bounded delay. Since we allow for the case where $R(0)=0$, it is possible that the delay becomes arbitrarily large if the phytoplankton population approaches zero. However, for any given state, the delay is bounded for all states within its neighbourhood. Results for differential equations with locally bounded delay can be found in \cite{walther}. However, our system is in the form that it my be transformed into a system with state-independent delay as in \cite{smith1991, smith1993}. We will therefore use this approach.

\section{Transformation to a State-Independent Delay Differential Equation}

It is possible to remove the state-dependence on the delay through a clever change in variable, as given in \cite{smith1991, smith1993}. The transformation is a state-dependent change of the time variable:
\begin{subequations}
\label{transform}
\begin{align}
    \hat{t} &= \int_0^t \frac{R(P(u))}{R^*}\,du, \label{time_transform}\\
    \hat{N}(\hat{t}) &= N(t), \label{N_transform}\\
    \hat{P}(\hat{t}) &= P(t), \label{P_transform}\\
    \hat{Z}(\hat{t}) &= Z(t), \label{Z_transform}
\end{align}
\end{subequations}
for $t>0$, where $R^*$ is a typical value of $R(P)$. With this transformation $t-\tau(s,P_t)$ becomes a constant, since for $\hat{t}\geq s/R^*$ we have
\begin{align*}
    \hat{t}-\frac{s}{R^*} = \int_0^t \frac{R(P(u))}{R^*}\,du-\int_{t-\tau(s,P_t)}^t \frac{R(P(u))}{R^*}\,du =\int_0^{t-\tau(s,P_t)} \frac{R(P(u))}{R^*}\,du,
\end{align*}
from the definition of $\tau$ in equation \eqref{tau_equation}. Then $\hat{P}(\hat{t}-s/R^*)=P(t-\tau(s,P_t))$ and $\hat{Z}(\hat{t}-s/R^*)=Z(t-\tau(s,P_t))$. That is, the delay no longer depends on the state. Furthermore, we can obtain an explicit equation for $\tau$. Define $\hat{\tau}:[0,m]\times \mathcal{C}[-m/R^*,0]\rightarrow [0,\infty)$ through
\begin{align}
\label{tauhat}
    \hat{\tau}(s,\hat{P}_{\hat{t}}) = \tau(s,P_t) = \int_{t-\tau(s,P_t)}^t \hspace{-1cm} \,du \hspace{0.4cm}= \int_{\hat{t}-s/R^*}^{\hat{t}} \frac{R^*}{R(\hat{P}(r))} \,dr,
\end{align}
where we have used the substitution
\begin{align*}
    r = \int_0^u \frac{R(P(v))}{R^*} \,dv.
\end{align*}
In a similar way, we can invert the time transformation:
\begin{align}
\label{time_invert}
    t = \int_0^t \,du = \int_0^{\hat{t}} \frac{R^*}{R(\hat{P}(r))} \,dr.
\end{align}

With this change of variables, the TDE model \eqref{model_TDE_2} is transformed to the following delay differential equation:
\begin{subequations}
\label{model_FDE}
\begin{align}
\frac{d\hat{N}(\hat{t})}{d\hat{t}} =& \frac{R^*}{R(\hat{P}(\hat{t}))} [-\mu \hat{P}(\hat{t})f(\hat{N}(\hat{t})) + \lambda \hat{P}(\hat{t})+ \delta \hat{Z}(\hat{t}) +(1-\gamma)g\hat{Z}(\hat{t})h(\hat{P}(\hat{t})) \nonumber\\
                &+ \delta_0(N_T-\hat{N}(\hat{t})-\hat{P}(\hat{t})-\hat{Z}(\hat{t}))], \label{model_FDE_N} \\
\frac{d\hat{P}(\hat{t})}{d\hat{t}} =& \frac{R^*}{R(\hat{P}(\hat{t}))} [\mu \hat{P}(\hat{t}) f(\hat{N}(\hat{t})) -\lambda \hat{P}(\hat{t}) - g \hat{Z}(\hat{t}) h(\hat{P}(\hat{t}))], \label{model_FDE_P}\\
\frac{d\hat{Z}(\hat{t})}{d\hat{t}} =& \gamma g e^{-\delta_0 \hat{\tau}(m,\hat{P}_{\hat{t}})} \frac{R^*}{R(\hat{P}(\hat{t}-T))}\hat{Z}(\hat{t}-T) h(\hat{P}(\hat{t}-T))-\frac{R^*}{R(\hat{P}(\hat{t}))} \delta \hat{Z}(\hat{t}), \label{model_FDE_Z}
\end{align}
\end{subequations}
for $\hat{t}\geq T$, where $T=m/R^*$ and $\hat{\tau}(m,\hat{P}_{\hat{t}})$ is defined by equation \eqref{tauhat}. We will refer to this system as the DDE model. The corresponding initial conditions are
\begin{align}
\label{model_FDE_IC}
    \hat{N}(T+\hat{t}) = \hat{\phi}_1(\hat{t}), \hspace{0.2cm} \hat{P}(T+\hat{t})=\hat{\phi}_2(\hat{t}),\hspace{0.2cm} \hat{Z}(T+\hat{t})=\hat{\phi}_3(\hat{t}),
\end{align}
for $ \hat{t}\in [-T,0]$, where $\hat{\phi}_1,\hat{\phi}_2,\hat{\phi}_3 \in \mathcal{C}[-T,0]$.  Again, solutions corresponding to the PDE1 model \eqref{model_PDE} should satisfy
\begin{align}
\label{conlaw_FDE_IC}
    N_T=\hat{\phi}_1(0)+\hat{\phi}_2(0)+\hat{\phi}_3(0)+\int_0^m e^{-\delta_0 \hat{\tau}(s,\hat{\phi}_2)}\frac{\gamma g \hat{\phi}_3(-\frac{s}{R^*})h(\hat{\phi}_2(-\frac{s}{R^*}))}{R(\hat{\phi}_2(-\frac{s}{R^*}))}\,ds.
\end{align}
For an initial condition that satisfies equation \eqref{conlaw_FDE_IC}, the solution to system \eqref{model_FDE} satisfies
\begin{align}
\label{conlaw_FDE}
    N_T=\hat{N}(\hat{t})+\hat{P}(\hat{t})+\hat{Z}(\hat{t})+\int_0^m e^{-\delta_0 \hat{\tau}(s,\hat{P}_{\hat{t}})}\frac{\gamma g \hat{Z}(\hat{t}-\frac{s}{R^*})h(\hat{P}(\hat{t}-\frac{s}{R^*}))}{R(\hat{P}(\hat{t}-\frac{s}{R^*}))}\,ds,
\end{align}
for all $\hat{t}\geq T$ for which the solution exists.

If we were to apply the time transformation in equation \eqref{time_transform} directly to the PDE2 model \eqref{model_PDE_2}, and then proceed with the method of characteristics, we would arrive at the DDE model \eqref{model_FDE}. In essence, the time transformation removes the dependency of the characteristic curves on $P(t)$, and instead places this dependence on the solutions along the curves. In \cite{smith1994} and \cite{smith1995}, he applies an analogous time transform to a system coupling a single adult population coupled with its corresponding maturity-structured juvenile population.

Since the delay in the DDE model \eqref{model_FDE} is bounded and state-independent, the standard theory for functional differential equations \cite{halebook} may be applied.

\section{Existence and Uniqueness of Solutions}

In this section, we will develop conditions for there to be a unique solution to the initial value problem for the TDE model \eqref{model_TDE_IC}-\eqref{model_TDE_2} that exists for all time. This will be done by matching it to a corresponding initial value problem for the DDE model \eqref{model_FDE}-\eqref{model_FDE_IC}. First, define the open set
\begin{align}
\label{phase_space}
    \hat{\Omega} = \{(\hat{\psi}_1,\hat{\psi}_2,\hat{\psi}_3)^T\in \mathcal{C}[-T,0] : R(\hat{\psi}_2(\theta)) \neq 0 \mbox{ for } \theta\in [-T,0] \}.
\end{align}
Standard results on functional differential equations with bounded delay (for example, see \cite{halebook}) give us that for $(\hat{\phi_1},\hat{\phi_2},\hat{\phi_3})^T\in \hat{\Omega}$, the initial value problem \eqref{model_FDE}-\eqref{model_FDE_IC} has a unique maximal solution.  We can then deduce the following result.

\begin{proposition}
\label{xhat_bounded}
    If $\hat{\phi}_1,\hat{\phi}_2,\hat{\phi}_3$ are positive functions that satisfy equation \eqref{conlaw_FDE_IC}, then the unique maximal solution to the initial value problem \eqref{model_FDE}-\eqref{model_FDE_IC}, denoted as $(\hat{N},\hat{P},\hat{Z})^T: [0,\hat{w})\rightarrow \mathbb{R}^3$, satisfies $\hat{N}(\hat{t}),\hat{P}(\hat{t}),\hat{Z}(\hat{t}) \in (0,N_T)$ for all $\hat{t} \in [0,\hat{w})$.
\end{proposition}
\begin{proof}
    Since $\hat{\phi}_1,\hat{\phi}_2,\hat{\phi}_3$ are positive functions, the result is clearly true for $\hat{t} \in [0,T]$.
    By the properties of $R$ in equation \eqref{Rprops} and the definition of the phase space, $\hat{\Omega}$ in equation \eqref{phase_space}, we find that $\hat{P}(\hat{t})>0$ for all $\hat{t} \in [0,\hat{w})$. Suppose there exists $\beta_1>T$ such that $\hat{Z}(\beta_1)=0$ and $\hat{Z}(\hat{t})>0$ for $\hat{t}<\beta_1$. Then, for $\hat{t}\in [T,\beta_1]$ we have that $\frac{d\hat{Z}(\hat{t})}{d\hat{t}}\geq-\delta \hat{Z}(\hat{t})$, which implies that $\hat{Z}(\hat{t})>0$ since $\hat{Z}(T)>0$. In particular, $\hat{Z}(\beta_1)>0$, which is a contradiction, so it must be true that $\hat{Z}(\hat{t})>0$ for all $\hat{t} \in [0,\hat{w})$. Equation \eqref{conlaw_FDE} then implies that $\hat{N}(\hat{t})+\hat{P}(\hat{t})+\hat{Z}(\hat{t})<N_T$, which implies that $\hat{N}$ is increasing at sufficiently small and positive values of $\hat{N}(\hat{t})$. Therefore, $\hat{N}(\hat{t})>0$ for all $\hat{t} \in [0,\hat{w})$. From equation \eqref{conlaw_FDE},  $\hat{N}(\hat{t}),\hat{P}(\hat{t}),\hat{Z}(\hat{t})>0$ implies also that $\hat{N}(\hat{t}),\hat{P}(\hat{t}),\hat{Z}(\hat{t})<N_T$. \eop
\end{proof}

It follows from this Proposition that all solutions are bounded. From the definition of the phase space, $\hat{\Omega}$ in equation \eqref{phase_space}, standard results for functional differential equations \cite{halebook} imply that each solution to the initial value problem \eqref{model_FDE}-\eqref{model_FDE_IC} either exists for all time or satisfies
\begin{align*}
    \lim_{\hat{t}\rightarrow \hat{w}^-} R(\hat{P}(\hat{t})) = 0.
\end{align*}
That is, solutions either exist for all time or approach $\partial \hat{\Omega}$.

\begin{proposition}
\label{bounded}
    If $(\phi_1,\phi_2,\phi_3)^T \in D_{N_T}$, then the solution to the initial value problem \eqref{model_TDE_IC}-\eqref{model_TDE_2} exists for all time and satisfies $N(t),P(t),Z(t)\in (0,N_T)$ for all $t\in [0,\infty)$.
\end{proposition}
\begin{proof}
    Consider the corresponding initial value problem \eqref{model_FDE}-\eqref{model_FDE_IC} with $\hat{\phi}_i(\hat{t})=\phi_i(t)$ for $i=1,2,3$ with $\hat{t}$ and $t$ related through equation \eqref{time_transform}. There exists a unique maximal solution $(\hat{N},\hat{P},\hat{Z})^T:[0,\hat{w})\rightarrow \mathbb{R}^3$, where $\hat{w}\in (T,\infty]$. Using the results in \cite{smith1991}, $(N,P,Z)^T:[0,w)\rightarrow \mathbb{R}^3$ defined by $(N(t),P(t),Z(t))^T=(\hat{N}(\hat{t}),\hat{P}(\hat{t}),\hat{Z}(\hat{t}))^T$ is the unique maximal solution to the initial value problem \eqref{model_TDE_IC}-\eqref{model_TDE_2} with $w = \lim_{\hat{t}\rightarrow \hat{w}^-}\int_0^{\hat{t}}R^*(R(\hat{P}(r)))^{-1}\,dr$.

    If $\hat{w} = \infty$, then
    \begin{align*}
       w= \lim_{\hat{t}\rightarrow \infty }\int_{0}^{\hat{t}} \frac{R^*}{R(\hat{P}(r))}\,dr = \infty
    \end{align*}
    since $R^*/R(\hat{P}(r))>R^*/R(N_T)>0$.

    If $\hat{w}<\infty$, then we have that $\lim_{\hat{t}\rightarrow \hat{w}^-} (\hat{N}(\hat{t}),\hat{P}(\hat{t}),\hat{Z}(\hat{t}))^T \in \partial \hat{\Omega}$ \cite{halebook}.  Since $(\hat{N}(\hat{t}),\hat{P}(\hat{t}),\hat{Z}(\hat{t}))^T$ is bounded by Proposition \ref{xhat_bounded}, in order for the solution to approach $\partial \hat{\Omega}$ we must have that
    \begin{align*}
    \lim_{\hat{t}\rightarrow \hat{w}^-} R(\hat{P}(\hat{t})) = 0.
    \end{align*}
    This implies that $\lim_{\hat{t}\rightarrow \hat{w}^-} \hat{P}(\hat{t}) = 0.$ From equation \eqref{model_FDE_P} and the properties of $h$ and $R$ in equations \eqref{h_props1} and \eqref{Rprops}, there exists a positive A such that
    \begin{align*}
        \frac{d}{d\hat{t}}R(\hat{P}(\hat{t})) > -A
    \end{align*}
    for $\hat{t} \in (\hat{w}-T,\hat{w})$. It then follows that $R(\hat{P}(\hat{t}))<A(\hat{w}-\hat{t})$ for $\hat{t} \in (\hat{w}-T,\hat{w})$, which implies that
    \begin{align*}
       w = \lim_{\hat{t}\rightarrow \hat{w}^- }\int_{0}^{\hat{t}} \frac{R^*}{R(\hat{P}(r))}\,dr \geq \lim_{\hat{t}\rightarrow \hat{w}^- }\int_{\hat{t}-T}^{\hat{t}} \frac{R^*}{A(\hat{w}-r)}\,dr = \infty.
    \end{align*}

    The last part, that $N(t),P(t),Z(t)\in (0,N_T)$ for all $t\in [0,\infty)$, follows directly from Proposition \ref{xhat_bounded} and the fact that $(N(t),P(t),Z(t))^T=(\hat{N}(\hat{t}),\hat{P}(\hat{t}),\hat{Z}(\hat{t}))^T$. \eop
\end{proof}

\section{Equilibrium Solutions}

To begin, note that $(N^*,P^*,Z^*)^T$ is an equilibrium solution of the DDE model \eqref{model_FDE} if and only if $(N^*,P^*,Z^*)^T$ with $R(P^*)\neq 0$ is an equilibrium solution of the TDE model \eqref{model_TDE_2}. For a fixed value of $N_T$, the system \eqref{model_FDE} has an equilibrium solution $(\hat{N}(\hat{t}),\hat{P}(\hat{t}),\hat{Z}(\hat{t}))^T=(N^*,P^*,Z^*)^T$ for $t\geq 0$ if $N^*,P^*,Z^*$ are constants that satisfy
\begin{subequations}
\label{newmodel_eq}
\begin{align}
\frac{R^*}{R(P^*)}[-\mu P^*f(N^*) + \lambda P^*+ \delta Z^* +(1-\gamma)gZ^*h(P^*) \nonumber\\
+ \delta_0 (N_T-N^*-P^*-Z^*)] &= 0,\label{newmodel_eq_N}\\
\frac{R^*}{R(P^*)}[\mu P^* f(N^*) -\lambda P^* - g Z^* h(P^*)] &= 0, \label{newmodel_eq_P}\\
\frac{R^*}{R(P^*)}[\gamma ge^{-\delta_0 m/R(P^*)} Z^* h(P^*)-\delta Z^*]&= 0 \label{newmodel_eq_Z}.
\end{align}
\end{subequations}
If $\delta_0=0$, then these three equations are redundant in the sense that if two are satisfied, then the third is also true. To avoid this problem, we use the conservation law \eqref{conlaw_TDE} in place of equation \eqref{newmodel_eq_N}. For $P^*\neq0$, equilibrium solutions should satisfy
\begin{subequations}
\label{newmodel_eq2}
\begin{align}
N^*+P^*+Z^*+\gamma g Z^* h(P^*)\frac{1-e^{-\delta_0 m/R(P^*)}}{\delta_0} -N_T &= 0,\label{newmodel_eq_N2}\\
\mu P^* f(N^*) -\lambda P^* - g Z^* h(P^*) &= 0, \label{newmodel_eq_P2}\\
\gamma ge^{-\delta_0 m/R(P^*)} Z^* h(P^*)-\delta Z^*&= 0 \label{newmodel_eq_Z2},
\end{align}
\end{subequations}
with $\frac{1-e^{-\delta_0 m/R(P^*)}}{\delta_0}=m/R(P^*)$ when $\delta_0=0$. If $\delta_0 \neq 0$  and $P^* \neq 0$, equations \eqref{newmodel_eq_N}-\eqref{newmodel_eq_Z} are satisfied if and only if equations \eqref{newmodel_eq_N2}-\eqref{newmodel_eq_Z2} are satisfied. However, when $\delta_0=0$, equations \eqref{newmodel_eq_N2}-\eqref{newmodel_eq_Z2} give the equilibrium solution corresponding to the value of total biomass, $N_T$, in which we are interested, whereas equations \eqref{newmodel_eq_N}-\eqref{newmodel_eq_Z} do not.

There are two types of equilibrium solutions that can exist: $E_1=(N_1^*,P_1^*,0)^T$, and $E_2 = (N_2^*,P_2^*,Z_2^*)^T$. We will say that an equilibrium point exists if all its components are non-negative, as these represent physical quantities. We will only consider $Z_2^*>0$ to distinguish between the two types.

Note that any equilibrium solution $(N^*,P^*,Z^*)$ of the DDE or TDE model corresponds to an equilibrium solution $(N^*,P^*,Z^*,\rho^*(s))$ of the $PDE1$ or $PDE2$ model, with
\begin{align*}
    \rho^*(s) = \frac{\gamma g Z^* h(P^*)}{R(P^*)}e^{\frac{-\delta_0s}{R(P^*)}}.
\end{align*}
The PDE1 and PDE2 models admit the solution $(N,P,Z,\rho(.,s))^T \equiv (N_T,0,0,0)^T$. We cannot define a corresponding solution $(N,P,Z)^T \equiv (N_T,0,0)^T$ to the TDE model \eqref{model_TDE_2} or the DDE model \eqref{model_FDE} since neither model is defined for $P\equiv 0$. However, it is possible for these models to asymptotically approach the point $(N_T,0,0)^T$. For this reason, we will define the point $e_0=(N_T,0,0)^T$, where the lower case $e$ is for emphasis that this is not an equilibrium solution, but a limit point.

\subsection{The Phytoplankton-only Equilibrium $E_1$}
Considering $E_1$, we see that $Z^*=0$ implies that \eqref{newmodel_eq_Z2} is satisfied. Then, for positive $P_1^*$, \eqref{newmodel_eq_P2} is satisfied if and only if $N_1^* = f^{-1}(\lambda/\mu)$. Consequently, \eqref{newmodel_eq_N2} is satisfied if and only if $P_1^* = N_T - N_1^*$. Therefore, $E_1$ exists if and only if
\begin{align*}
    N_T > N_{T1},
\end{align*}
where
\begin{align}
\label{NT1}
    N_{T1} = f^{-1}\left(\frac{\lambda}{\mu}\right).
\end{align}
Note that as $N_T$ increases, $P_1^*$ increases linearly while $N_1^*$ remains fixed. That is, while $E_1$ is a stable equilibrium, increasing the total biomass in the ecosystem increases the phytoplankton population while the dissolved nutrient remains fixed.

\subsection{The Phytoplankton-Zooplankton Equilibrium $E_2$}
For the equilibrium $E_2$, we have that $Z_2^*>0$, so it is required that
\begin{align*}
\gamma ge^{-\delta_0 m/R(P_2^*)} h(P_2^*)-\delta = 0 ,
\end{align*}
in order for \eqref{newmodel_eq_Z2} to be satisfied. This is true if and only if
\begin{align}
\label{m_Pstar}
    m = \frac{R(P_2^*)}{\delta_0}\ln\left(\frac{\gamma g h(P_2^*)}{\delta}\right).
\end{align}
Given $m$, $P_2^*$ is defined implicitly through this equation.  Under the assumption that $R$ and $h$ are both increasing and saturating functions, we can see that $P_2^*$ increases with $m$, $P_2^* = h^{-1}(\delta/\gamma g)$ when $m=0$, and $P_2^* \rightarrow \infty$ as $m \nearrow R_{\infty}\ln(\gamma g/\delta)/\delta_0$. There is no solution for $m$ larger than this value. However, if $m \in [0,R_{\infty}\ln(\gamma g/\delta)/\delta_0)$, then there is a unique positive value $P_2^*$ that satisfies \eqref{m_Pstar}.

Assuming that $m \in [0,R_{\infty}\ln(\gamma g/\delta)/\delta_0)$, we find that \eqref{newmodel_eq_P2} is satisfied if and only if
\begin{align*}
    Z_2^* = (\mu f(N_2^*)-\lambda)\frac{P_2^*}{gh(P_2^*)}.
\end{align*}
For $Z_2^*>0$ , we require $N_2^*>N_{T1}$. Then \eqref{newmodel_eq_N2} is satisfied if $N_2^*$ satisfies
\begin{align*}
    N_2^*+P_2^*+(\mu f(N_2^*)-\lambda)\frac{P_2^*}{gh(P_2^*)}\left(1+\gamma g h(P_2^*)\frac{1-e^{-\delta_0 m/R(P_2^*)}}{\delta_0}\right) -N_T &= 0,
\end{align*}
with $P_2^*$ fixed and given by \eqref{m_Pstar}. We can see that $N_2^*$ increases with increasing $N_T$ and that $N_2^*=N_{T1}$ when $N_T=N_{T1}+P_2^*$. Therefore, $E_2$ exists and is unique if and only if $m \in [0,R_{\infty}\ln(\gamma g/\delta)/\delta_0)$ and $N_T>N_{T2}$ where
\begin{align*}
    N_{T2} = f^{-1}\left(\frac{\lambda}{\mu}\right) + P_2^*,
\end{align*}
with $P_2^*$ defined through equation \eqref{m_Pstar}. Since $N_{T2}>N_{T1}$, the existence of $E_2$ implies the existence of $E_1$.  As $N_T$ is increased, $P_2^*$ is fixed while $N_2^*$ and $Z_2^*$ increase.

\begin{figure}
\begin{center}
\includegraphics[width=0.45 \textwidth]{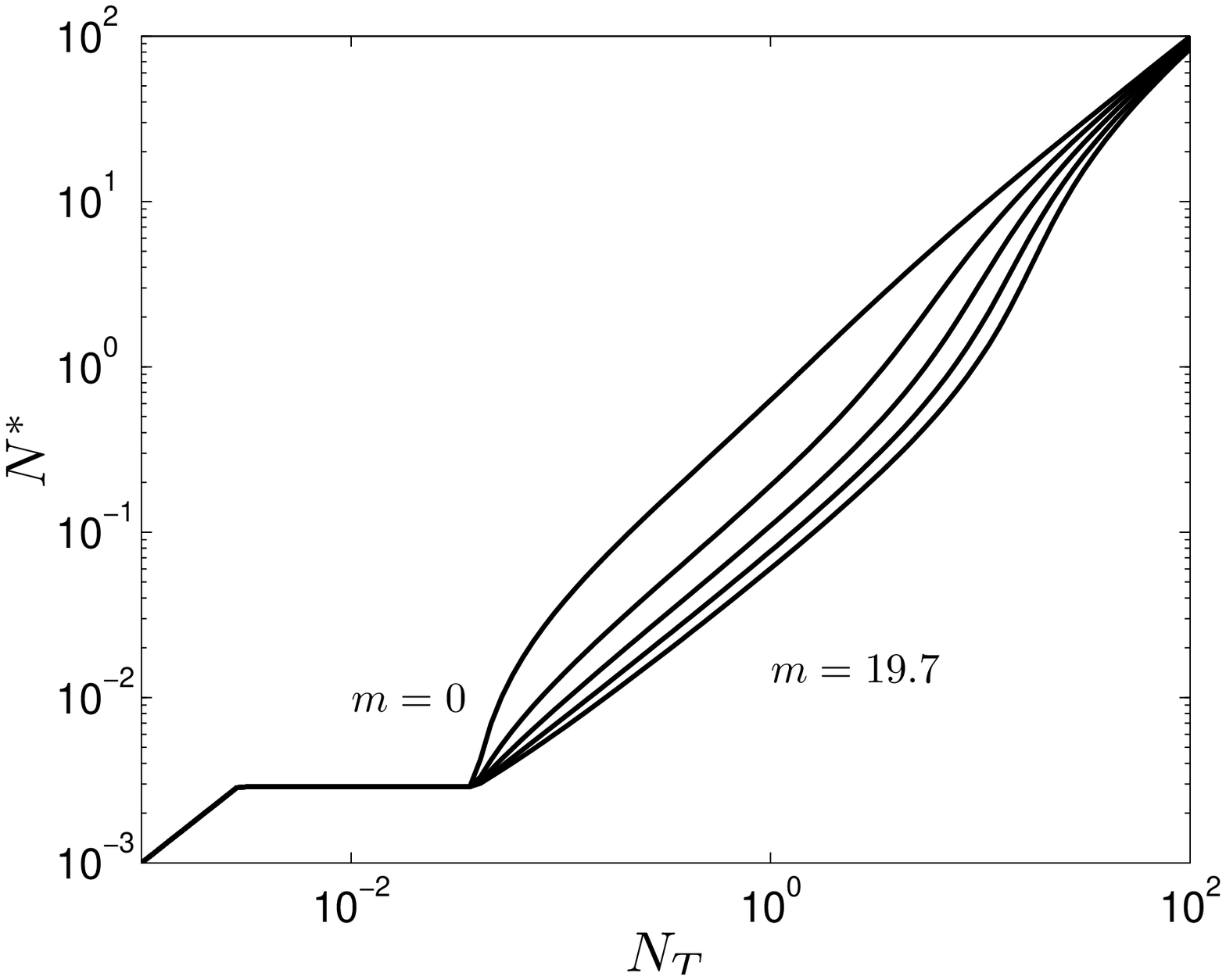}
\includegraphics[width=0.45 \textwidth]{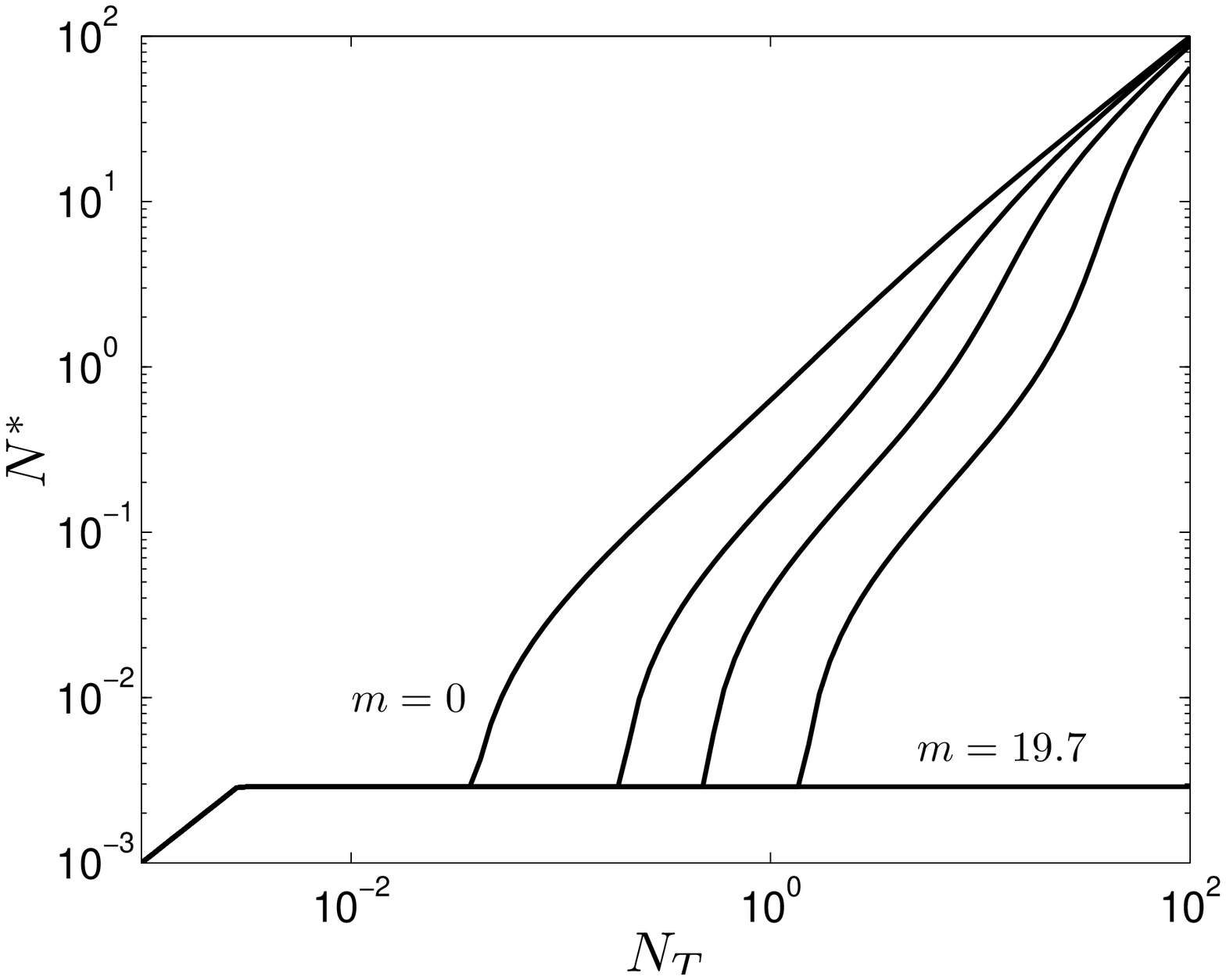}\\
\includegraphics[width=0.45 \textwidth]{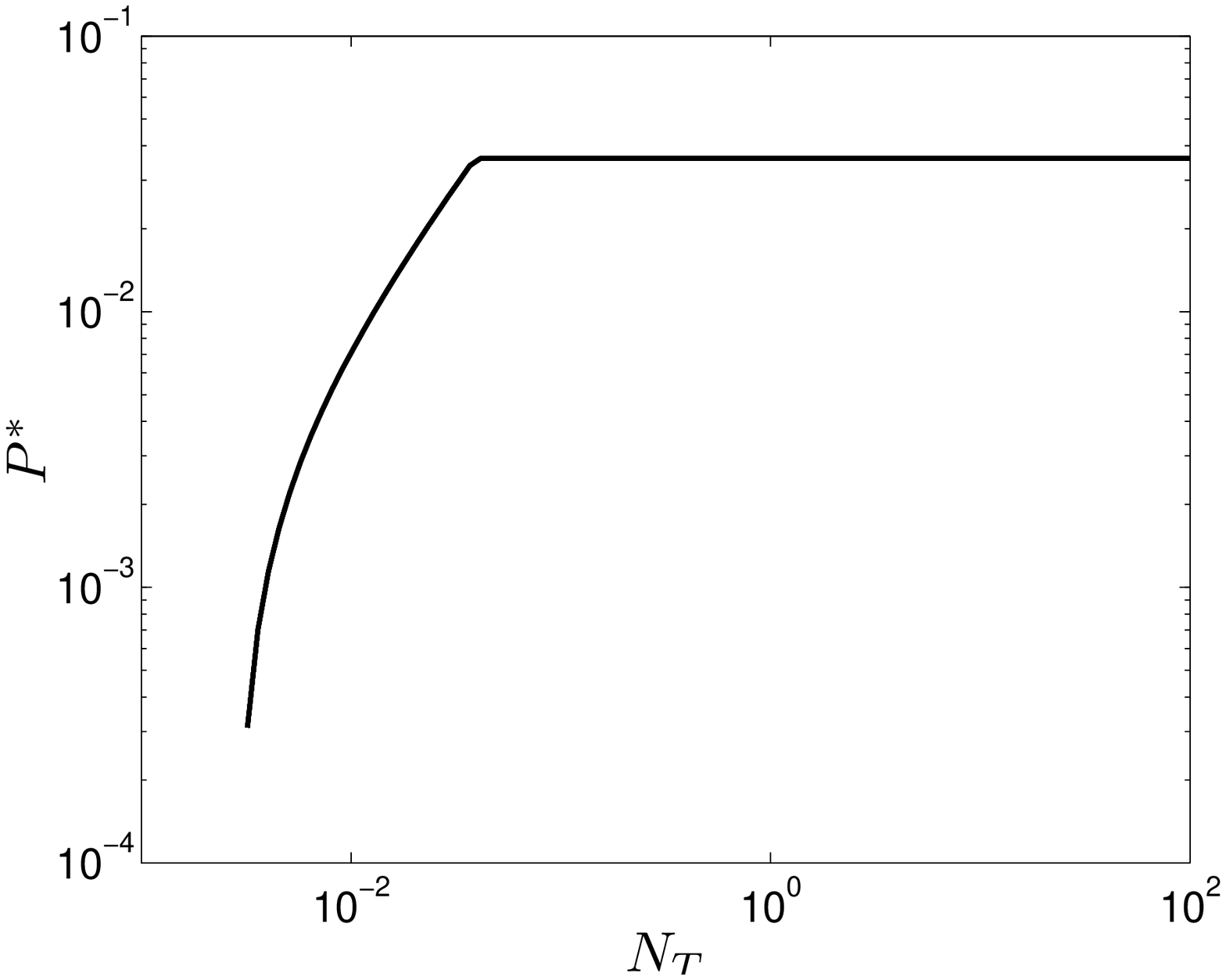}
\includegraphics[width=0.45 \textwidth]{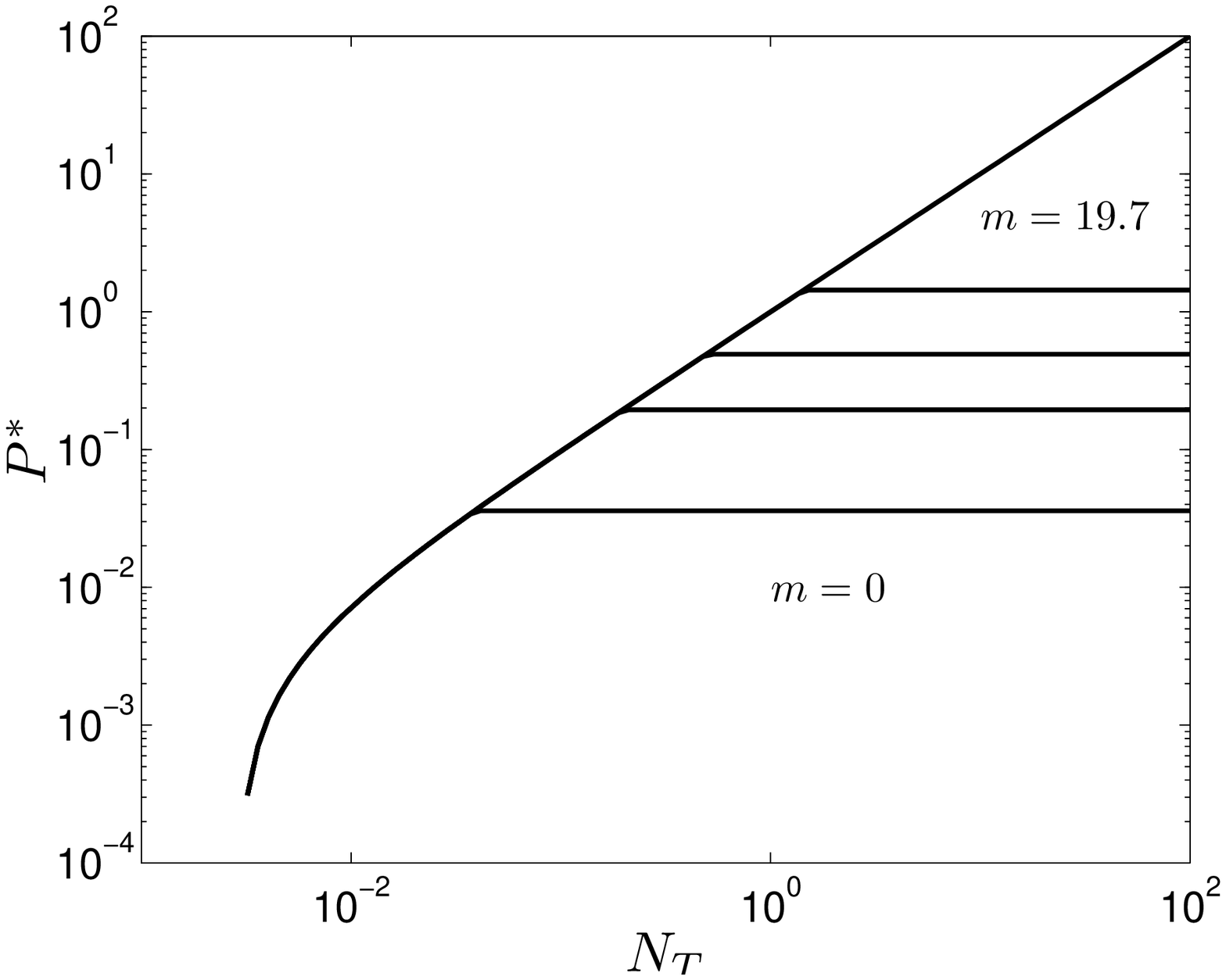}\\
\includegraphics[width=0.45 \textwidth]{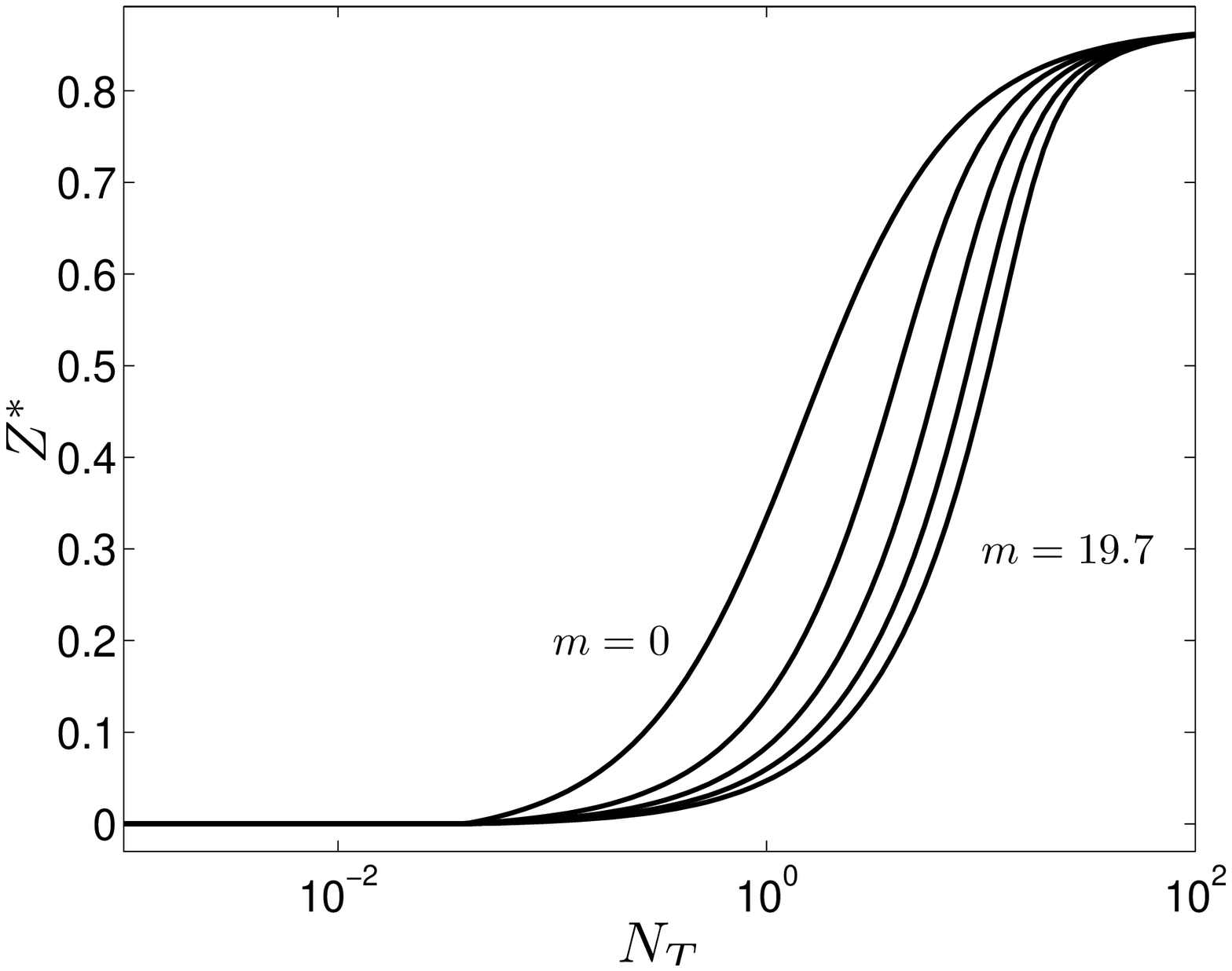}
\includegraphics[width=0.45 \textwidth]{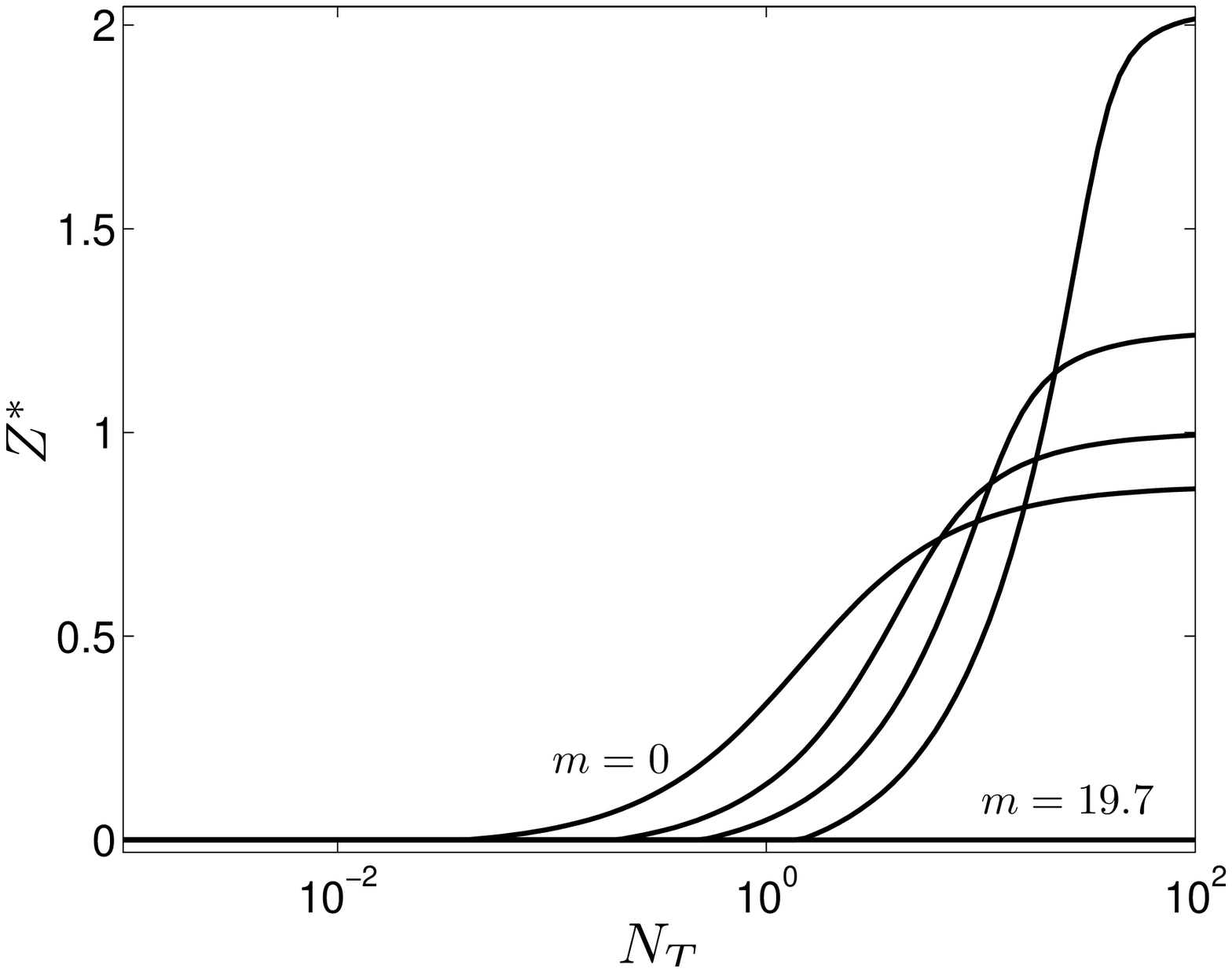}
\end{center}
\caption{Equilibrium solutions as a function of total biomass for $m=0,5,10,15,19.7$ and $R(P)=P/(P+0.159)$. The left three are for $\delta_0=0$ and the right three are for $\delta_0=\delta$. If $E_2$ exists, it is the solution plotted. If $E_2$ does not exist, but $E_1$ does, then $E_1$ is plotted. If neither exist, the limit point $(N_T,0,0)^T$ is plotted}
\label{equilsolutions}
\end{figure}

Figure \ref{equilsolutions} shows how the dominant equilibrium solutions change as the total biomass increases. We use parameter values given in Table \ref{params}, which were taken from \cite{poulin}. We consider the case where $R(P)=P/(P+0.159)$ and various values of $m$. The top three plots are for $\delta_0=0$ and the bottom three are for $\delta_0 = \delta$. We see that $N^*=N_T$ and $P^*=Z^*=0$ for $N_T < N_{T1}$ ($e_0$ is plotted). Then for $N_{T1}<N_T<N_{T2}$ ($E_1$ is plotted), $P^*$ increases linearly with $N_T$ while $N^*$ is fixed at $N_{T1}$ and $Z^*$ is fixed at 0. Then for $N_T>N_{T2}$ ($E_2$ is plotted), $P^*$ is fixed and $N^*$ and $Z^*$ increase with $N_T$, though $Z^*$ saturates and $N^*$ increases indefinitely. Note that $N_{T2}$ depends on $m$ when $\delta_0 \neq 0$, but not when $\delta_0=0$. Also, $P^*$ is independent of $m$ when $\delta_0=0$, but changes for $N_T>N_{T2}$ when $\delta_0 \neq 0$. The value to which $Z^*$ saturates does not change with $m$ when $\delta_0 =0$, but increases with $m$ when $\delta_0 \neq 0$.

\begin{table}[t]
\footnotesize
\begin{center}
\begin{tabular}{|c|c|}
\hline
Parameter & Value \\
\hline
$\mu$ & $5.9$ day$^{-1}$\\
$\lambda$ & $0.017$ day$^{-1}$\\
$g$ & $7$ day$^{-1}$\\
$\gamma$ & $0.7$\\
$\delta$ & $0.17$ day$^{-1}$\\
$f(N)$ & $\frac{N}{N+k} $ \\
$h(P)$ & $\frac{P}{P+K} $ \\
$k$ & $1.0 \mu\mbox{M} $ \\
$K$ & $1.0 \mu\mbox{M} $ \\
\hline
\end{tabular}
\caption{Parameter values used for all computations}
\end{center}
\label{params}
\end{table}


\section{Linearization and Characteristic Equation}

In this section, we will show that a general threshold delay equation and its corresponding delay differential equation, as in \cite{smith1991}, have identical linearizations about an equilibrium solution.

Consider the threshold delay equation in the general form
\begin{align}
\label{general_TDE}
    \frac{dx(t)}{dt} = f(x(t),x(t-\tau(x_t)),\tau(x_t)), \hspace{0.5cm} \int_{-\tau(x_t)}^0 K(x_t(u))\,du=m,
\end{align}
and the corresponding functional differential equation
\begin{align}
\label{general_FDE}
    \frac{d\hat{x}(\hat{t})}{d\hat{t}} = \frac{K^*}{K(\hat{x}(\hat{t}))}f(\hat{x}(\hat{t}),\hat{x}(\hat{t}-m/K^*),\hat{\tau}(\hat{x}_{\hat{t}})), \hspace{0.5cm} \hat{\tau}(\hat{x}_{\hat{t}})=\int_{-m/K^*}^0\frac{K^*}{K(\hat{x}_{\hat{t}}(u))} \,du,
\end{align}
where $f$ and $K$ are differentiable, and $K^*$ is a typical value of $K$. Results on the linearization of state-dependent delay differential equations with bounded delay can be found in \cite{cooke1996}. These results apply to systems with bounded delay, which do not apply necessarily apply directly to the system here. Nevertheless, we will show that applying the linearization procedure for a state-dependent delay system with bounded delay to system \eqref{general_TDE} yields the same linearization as linearizing system \eqref{general_FDE} in the standard way for fixed delay systems.

Let $x^*$ be an equilibrium solution of equations \eqref{general_TDE} and \eqref{general_FDE}. That is, a constant function such that $f(x^*,x^*,\tau^*)=0$ where $\tau^* = m/K(x^*)=\tau(x^*)=\hat{\tau}(x^*)$. Let $y(t)=x(t)-x^*$. By setting the state-dependent delay to its equilibrium value, and linearizing equation \eqref{general_TDE} as in \cite{cooke1996}, we obtain
\begin{align}
\label{general_TDE_lin}
    \frac{dy(t)}{dt} =& D_1f(x^*,x^*,\tau^*) y(t) +D_2f(x^*,x^*,\tau^*) y(t-\tau^*) \nonumber\\
    &+D_3f(x^*,x^*,\tau^*) D\tau(x^*) y_t
\end{align}
where $D_if$ is the derivative of $f$ with respect to its $i$th argument and
\begin{align}
\label{tau_diff}
    D\tau(x^*)y_t = -\frac{DK(x^*)}{K(x^*)}\int_{-\tau^*}^0 y_t(u) \,du.
\end{align}

Let $\hat{y}(\hat{t})=\hat{x}(\hat{t})-x^*$.  Then linearizing equation \eqref{general_FDE} in the standard way yields
\begin{align}
\label{general_FDE_lin}
    \frac{d\hat{y}(\hat{t})}{d\hat{t}} = &\frac{K^*}{K(x^*)}[D_1f(x^*,x^*,\tau^*) \hat{y}(\hat{t}) +D_2f(x^*,x^*,\tau^*) \hat{y}(\hat{t}-m/K^*)\nonumber \\
    &+D_3f(x^*,x^*,\tau^*) D\hat{\tau}(x^*) \hat{y}_{\hat{t}}]
\end{align}
where
\begin{align}
\label{tauhat_diff}
    D\hat{\tau}(x^*)\hat{y}_{\hat{t}} = -\frac{K^*DK(x^*)}{K(x^*)^2}\int_{-m/K^*}^0 \hat{y}_{\hat{t}}(u) \,du.
\end{align}
Without any loss of generality, we can set $K^*=K(x^*)$. Indeed, this would be the natural choice of $K^*$ if we were interested in solutions near a given equilibrium. Then it can be see that equations \eqref{general_TDE_lin} and \eqref{general_FDE_lin} are identical.

We can check that \eqref{tau_diff} is correct by  verifying that
\begin{align}
\label{tau_limit}
    \lim_{h\rightarrow 0}\frac{\left|\tau(x^*+h)-\tau(x^*)+\frac{DK(x^*)}{K(x^*)} \int_{-\tau^*}^{0} h(u) \,du\right|}{||h||} = 0,
\end{align}
where $h \in \mathcal{C}[-r,0]$ for some $r>\tau(\phi)$ for all $\phi$ in a neighbourhood of $x^*$. Here, $||.||$ is the usual sup norm. We proceed as follows:
\begin{align*}
    & \left|\tau(x^*+h)-\tau(x^*)+\frac{DK(x^*)}{K(x^*)} \int_{-\tau^*}^{0} h(u) \,du\right| \\
   =& \frac{1}{K(x^*)} \left|\int_{-\tau(x^*+h)}^{0}K(x^*)\,du-\int_{-\tau^*}^{0}K(x^*)\,du + DK(x^*) \int_{-\tau^*}^{0} h(u) \,du\right|.
\end{align*}
Under the assumption that $K\in \mathcal{C}^1$, it is true that $K(x^*) = K(x^*+x) - DK(x^*)x + G(x)$, with $G$ satisfying
\begin{align*}
    \lim_{x\rightarrow 0} \frac{|G(x)|}{||x||_E} = 0,
\end{align*}
where $||.||_E$ is the Euclidean norm on $\mathbb{R}^n$. Then
\begin{align*}
     &\int_{-\tau(x^*+h)}^{0}K(x^*)\,du \\
     =& \int_{-\tau(x^*+h)}^{0}K(x^*+h(u))\,du-DK(x^*)\int_{-\tau(x^*+h)}^{0}h(u)\,du \\
     &+\int_{-\tau(x^*+h)}^{0}G(h(u))\,du.
\end{align*}
Then by the definition of $\tau$:
\begin{align*}
     \int_{-\tau(x^*+h)}^{0}K(x^*+h(u))\,du = \int_{-\tau^*}^{0}K(x^*)\,du = m.
\end{align*}
So it can then be seen that
\begin{align*}
    &\lim_{h\rightarrow 0} \frac{\left|\tau(x^*+h)-\tau(x^*)+\frac{DK(x^*)}{K(x^*)} \int_{-\tau(x^*)}^{0} h(u) \,du\right|}{||h||} \\
    =&\lim_{h\rightarrow 0} \frac{\left|DK(x^*)\int_{-\tau^*}^{-\tau(x^*+h)}h(u)\,du+\int_{-\tau(x^*+h)}^{0}G(h(u))\,du\right|}{K(x^*)||h||}\\
    =&0,
\end{align*}
from the properties of $G$.

To verify equation \eqref{tauhat_diff}, we check that
\begin{align}
\label{tau_limit_2}
    \lim_{h\rightarrow 0}\frac{\left|\hat{\tau}(x^*+h)-\hat{\tau}(x^*)+\frac{K^*DK(x^*)}{K(x^*)^2} \int_{-m/K^*}^{0} h(u) \,du\right|}{||h||} = 0.
\end{align}
This follows from the definition of $\hat{\tau}$ in equation \eqref{general_FDE} and the fact that for $x$ near $x^*$
\begin{align}
    \frac{1}{K(x^*+x)} = \frac{1}{K(x^*)}-\frac{DK(x^*)}{K(x^*)^2}x+\hat{G}(x),
\end{align}
for some function $\hat{G}$ that satisfies
\begin{align*}
    \lim_{x\rightarrow 0} \frac{|\hat{G}(x)|}{||x||_E} = 0.
\end{align*}

In particular, the linearized system for our model (\eqref{model_TDE_2} or \eqref{model_FDE}) about an equilibrium solution $(N^*,P^*,Z^*)^T$ with $R^*$ taken to be $R(P^*)$ is
\begin{align}
\label{lin_system}
    \frac{dy(t)}{dt} = A_1 y(t) + A_2 y(t-T) + A_3 \int_{-T}^0y(t+u) \,du,
\end{align}
where
\begin{align*}
    A_1 = \left( \begin{array}{ccc}
    -\mu P^* a -\delta_0 \hspace{0.15cm}& -\mu c +\lambda +(1-\gamma)gZ^*b-\delta_0 \hspace{0.15cm}& \delta-\delta_0 +(1-\gamma)gd \\
    \mu P^* a \hspace{0.15cm}& \mu c-\lambda-gZ^*b \hspace{0.15cm}& -gd \\
    0 \hspace{0.15cm}& e^{-\delta_0 T}\gamma g Z^*d\frac{R'(P^*)}{R(P^*)} \hspace{0.15cm}& -\delta
    \end{array} \right),
\end{align*}
\begin{align*}
    A_2 = \left( \begin{array}{ccc}
    0 & 0 & 0 \\
    0 & 0 & 0 \\
    0 & e^{-\delta_0 T}\gamma g Z^*(b-\frac{R'(P^*)}{R(P^*)}d) & e^{-\delta_0 T} \gamma g d
    \end{array} \right),
\end{align*}
\begin{align*}
    A_3 = \left( \begin{array}{ccc}
    0 & 0 & 0 \\
    0 & 0 & 0 \\
    0 & \delta_0 e^{-\delta_0 T} \gamma g Z^* d \frac{R'(P^*)}{R(P^*)}& 0
    \end{array} \right),
\end{align*}
and $a=f'(N^*), b=h'(P^*), c=f(N^*)$, $d=h(P^*)$, and $T=m/R(P^*)=\tau(m,P^*)=\hat{\tau}(m,P^*)$.

Substituting $y(t) = v e^{st}$ into \eqref{lin_system}, with $v\in \mathbb{R}^3$, we can obtain the characteristic equation:
\begin{align}
\label{characteristic_equation}
    \det\left(s I - A_1-A_2 e^{-sT} + A_3\frac{(1-e^{-sT})}{s}\right) = 0.
\end{align}
If all values of $s$ that satisfy equation \eqref{characteristic_equation} have a negative real part, then $(N^*,P^*,Z^*)^T$ is an asymptotically stable equilibrium solution of system \eqref{model_FDE}  \cite{halebook}. Due to the equivalence of systems \eqref{model_TDE_2} and \eqref{model_FDE}, they are asymptotically stable under identical conditions. Thus, we can use the characteristic equation \eqref{characteristic_equation} to obtain asymptotic stability in system \eqref{model_TDE_2}. Therefore, if all values of $s$ that satisfy equation \eqref{characteristic_equation} have a negative real part, then $(N^*,P^*,Z^*)^T$ is an asymptotically stable equilibrium solution of system \eqref{model_TDE_2}.


\section{Stability and Extinction}

\subsection{Conditions for the Extinction of Plankton}

It is possible for the solutions to the TDE model \eqref{model_TDE_2} to approach the point $e_0=(N_T,0,0)^T$. To show this, we will use the following Lemma.
\begin{lemma}
\label{extinct}
    If $(N(t),P(t),Z(t))^T$ is a solution to system \eqref{model_TDE_2} such that \\$\lim_{t \rightarrow \infty} P(t) = 0$, then $\lim_{t \rightarrow \infty} (N(t),P(t),Z(t))^T=(N_T,0,0)^T$.
\end{lemma}
\begin{proof}
    Assume that $\lim_{t \rightarrow \infty} P(t) = 0$. The term $e^{-\delta_0 \tau(m,P_t)}$ is clearly bounded. Also, by Proposition \ref{bounded}, $Z(t-\tau(m,P_t))<N_T$, so this term is also bounded. From the properties of $h$ and $R$ in equations \eqref{h_props1} and \eqref{Rprops}, there is a function, $w$, such that $\lim_{t \rightarrow \infty}w(t) =0$ and $\frac{dZ(t)}{dt} < w(t) -\delta Z(t)$ for $t\geq t_0$. Integrating this inequality, we obtain
    \begin{align*}
        Z(t) < e^{-\delta (t-t_0)}Z(t_0) + e^{-\delta t} \int_{t_0}^t e^{\delta u} w(u) \,du.
    \end{align*}
    In the case where $\int_{t_0}^{\infty} e^{\delta u} w(u) \,du < \infty$ it is obvious the final term approaches zero as $t \rightarrow \infty$. Otherwise, we can apply l'H\^opital's rule and get the same result. Since $Z(t)>0$ by Proposition \ref{bounded}, we have that $\lim_{t \rightarrow \infty}Z(t) = 0$.

    The result then follows from the conservation law \eqref{conlaw_TDE} if it is true that
     \begin{align*}
        \lim_{t \rightarrow \infty} \int_0^m e^{-\delta_0 \tau(s,P_t)}\frac{\gamma g Z(t-\tau(s,P_t))h(P(t-\tau(s,P_t)))}{R(P(t-\tau(s,P_t)))} \,ds=0.
    \end{align*}
    If $\lim_{t\rightarrow \infty}\tau(s,P_t)<\infty$, then $ Z(t-\tau(s,P_t)) \rightarrow 0$ and the result immediately follows. If $\lim_{t\rightarrow \infty}\tau(s,P_t)=\infty$ then $e^{-\delta_0 \tau(s,P_t)}\rightarrow 0$ and the result follows. \eop
\end{proof}

This also implies that if $\hat{P}(\hat{t})$ in the solution to the DDE model \eqref{model_FDE} reaches zero in finite time, then $\hat{Z}(\hat{t})$ reaches zero at the same finite moment in time. This is because the transformation maps $t=\infty$ to $\hat{t}=\hat{w}<\infty$.

The following Proposition shows conditions for the extinction of plankton.
\begin{proposition}
    If  $0<N_T<N_{T1}$, where $N_{T1}$ is given in \eqref{NT1}, and $(\phi_1,\phi_2,\phi_3)^T\in D_{N_T}$, then the solution to the initial value problem \eqref{model_TDE_IC}-\eqref{model_TDE_2} asymptotically approaches $(N_T,0,0)^T$.
\end{proposition}
\begin{proof}
    From \eqref{model_TDE2_P}, we have that $\frac{dP}{dt} < \mu f(N_T) P(t) -\lambda P(t)$. Then $N_T<N_{T1}$ implies that $\mu f(N_T)-\lambda<0$, which implies that $P(t)< P(t_0) e^{(\mu f(N_T)-\lambda)( t-t_0)}$. Since $P(t)>0$ from Proposition \ref{bounded}, $\lim_{t\rightarrow \infty} P(t) = 0$. The result then follows from Lemma \ref{extinct}. \eop
\end{proof}



\subsection{Global Stability of $E_1$}

Here we will show conditions for the global stability of $E_1$. Consider the following Lemma.

\begin{lemma}
\label{Pbeta}
    Let $(\phi_1,\phi_2,\phi_3)^T \in D_{N_T}$ and denote $(N,P,Z)^T: [0,\infty)\rightarrow \mathbb{R}^3$ as the solution to the initial value problem \eqref{model_TDE_IC}-\eqref{model_TDE_2}. If $N_{T}>N_{T1}$, then for any constant $\beta>N_T-N_{T1}$, there exists a $t_1\geq t_0$ such that $P(t)<\beta$ for all $t\geq t_1$.
\end{lemma}
\begin{proof}
    Fix $\beta < N_T-N_{T1}$. From the conservation law \eqref{conlaw_TDE} and Proposition \ref{bounded}, $N(t)<N_T-P(t)$. Then for $P(t)\geq \beta$ we have
    \begin{align*}
        \frac{dP(t)}{dt} &< [\mu f(N(t))-\lambda]P(t), \\
        &< [\mu f(N_T-P(t))-\lambda]P(t), \\
        &< [\mu f(N_T-\beta)-\lambda]P(t).
    \end{align*}
    Since $\mu f(N_T-\beta)-\lambda<\mu f(N_{T1})-\lambda=0$, we have
    \begin{align*}
        \frac{dP(t)}{dt} &< [\mu f(N_T-\beta)-\lambda]\beta
    \end{align*}
    for $P(t)\geq \beta$. Therefore there exists $t_1\geq t_0$ such that $P(t)<\beta$ for all $t\geq t_1$. \eop
\end{proof}

\begin{proposition}
\label{E1}
    Let $(\phi_1,\phi_2,\phi_3)^T \in D_{N_T}$ and denote $(N,P,Z)^T: [0,\infty)\rightarrow \mathbb{R}^3$ as the solution to the initial value problem \eqref{model_TDE_IC}-\eqref{model_TDE_2}. If $N_{T1}<N_{T}<N_{T2}$, then  $\lim_{t\rightarrow \infty} (N(t),P(t),Z(t))^T = (N_1^*,P_1^*,0)^T$.
\end{proposition}
\begin{proof}
    Note that $\mu f(N_1^*)-\lambda=0$ and $P_1^*=N_T-N_1^*<P_2^*$ for $N_T<N_{T2}$. Consider the function
    \begin{align*}
        V_1(t) = \int_{N_1^*}^{N(t)}[\mu f(x)-\lambda ] \,dx + \delta_0 \int_{N_T-N_1^*}^{P(t)}\frac{x+N_1^*-N_T}{x} \,dx.
    \end{align*}
    We have that $V_1(t)\geq 0$ for $P(t)>0$. Taking the derivative along solutions to the TDE model \eqref{model_TDE_2}, we obtain
    \begin{align*}
        \frac{d}{dt}V_1(t) =& [\mu f(N(t))-\lambda][-\mu P(t)f(N(t)) +\lambda P(t) + \delta Z(t) + (1-\gamma)gZ(t)h(P(t)) \\
    &+ \delta_0(N_T-N(t)-P(t)-Z(t))] \\ &+ \delta_0\frac{P(t)+N_1^*-N_T}{P(t)} [\mu P(t)f(N(t)) - \lambda P(t) - g Z(t) h(P(t))],\\
    \leq& -P(t)[\mu f(N(t))-\lambda]^2 -\delta_0[N(t)-N_1^*][\mu f(N(t))-\lambda]+M Z(t),
    \end{align*}
    for some positive constant $M$. This $M$ exists because $N(t)$ and $P(t)$ are bounded and from the properties of $h$ in \eqref{h_props1}. Then consider the function
    \begin{align*}
        V_2(t) = Z(t) + \int_{t-\tau(m,P_t)}^{t} \gamma g e^{-\delta_0 \tau(m,\beta)} Z(u)h(P(u)) \,du.
    \end{align*}
    We have that $V_2(t)\geq 0$. Choose $\beta \in (N_T-N_{T1},P_2^*)$. Then by Lemma \ref{Pbeta}, there exists $t_1\geq t_0$ such that $P(t)<\beta$ for $t\geq t_1$. This implies that
    \begin{align*}
        \gamma ge^{-\delta_0 \tau(m,P_t)} h(P(t))<\gamma ge^{-\delta_0 \tau(m,\beta)} h(\beta)<\gamma ge^{-\delta_0 \tau(m,P_2^*)} h(P_2^*)=\delta
    \end{align*}
    for $t\geq t_1$. Define the positive constant $\alpha = \delta -\gamma ge^{-\delta_0 \tau(m,\beta)} h(\beta)$. Differentiating, we have that
    \begin{align*}
        \frac{d}{dt}V_2(t) =& R(P(t)) e^{-\delta_0 \tau(m,P_t)}\frac{\gamma g Z(t-\tau(m,P_t))h(P(t-\tau(m,P_t)))}{R(P(t-\tau(m,P_t)))} - \delta Z(t)\\
        &+ \gamma g e^{-\delta_0 \tau(m,\beta)} Z(t)h(P(t))\\
        &-\gamma g e^{-\delta_0 \tau(m,\beta)} Z(t-\tau(m,P_t)h(P(t-\tau(m,P_t))\left(1-\frac{\partial \tau}{\partial t}\right), \\
        \leq& -\alpha Z(t),
    \end{align*}
    where we have used that
    \begin{align*}
        1-\frac{\partial \tau}{\partial t} = \frac{R(P(t))}{R(P(t-\tau(m,P_t))},
    \end{align*}
    $h(P(t))< h(\beta)$, and $e^{-\delta_0 \tau(m,P_t)} < e^{-\delta_0 \tau(m,\beta)}$ for $t\geq t_1$. Setting $V(t)=V_1(t)+\frac{M+M_1}{\alpha}V_2(t)$ where $M_1>0$, we have that $V(t)\geq 0$ and that
    \begin{align*}
        \frac{d}{dt}V(t) \leq -\delta_0[N(t)-N_1^*][\mu f(N(t))-\lambda]-M_1Z(t),
    \end{align*}
    for $t\geq t_1$. Integrating, we obtain that
    \begin{align}
    \label{V}
        V(t)+\int_{t_1}^{t}[\delta_0(N(u)-N_1^*)(\mu f(N(u))-\lambda)+M_1Z(u)]\,du \leq V(t_1).
    \end{align}
    Since $V(t)$ and the above integrand are nonnegative, and since the inequality \eqref{V} must be true for all time, it follows that
    \begin{align*}
        \int_{t_0}^{\infty}[\delta_0(N(u)-N_1^*)(\mu f(N(u))-\lambda)+M_1Z(u)]\,du<\infty.
    \end{align*}
    By Barb\u{a}lat's Lemma \cite{gopalsamy}, it must be true that $N(t)\rightarrow N_1^*$ and $Z(t)\rightarrow 0$. By the conservation law \eqref{conlaw_TDE} it then follows that $P(t)\rightarrow P_1^*$. \eop
\end{proof}

The attractivity of $E_1$ has thus been shown. However, Proposition \ref{E1} does not address its stability. Let $(N,P,Z)^T : [0,\infty)\rightarrow \mathbb{R}^3$ denote the solution to the initial value problem \eqref{model_TDE_IC} - \eqref{model_TDE_2}. We will consider an equilibrium point $E$ to be asymptotically stable if for any $\varepsilon>0$, there exists a neighbourhood $N_{\varepsilon}$ around $E$  such that if $(\phi_1,\phi_2,\phi_3) \in D_{N_T} \cap N_{\varepsilon}$ then $||(N(t),P(t),Z(t))^T -E|| < \varepsilon$ for $t\geq t_0$ and $\lim_{t\rightarrow \infty}(N(t),P(t),Z(t))^T  = E$.

We will use the linearization of the TDE model \eqref{model_TDE_2} to prove the following Proposition.
\begin{proposition}
    If $N_{T1}<N_T<N_{T2}$, then $E_1$ is asymptotically stable.
\end{proposition}
\begin{proof}
    The condition $N_{T1}<N_T$ ensures that $E_1=(N_1^*,P_1^*,0)^T$ exists. The characteristic equation of the corresponding linear system is
    \begin{align*}
    &\det \left( \begin{array}{ccc}
    s+\mu P_1^* a +\delta_0 \hspace{0.2cm} &  \delta_0\hspace{0.2cm}& -\delta +\delta_0 -(1-\gamma)gd\\
    -\mu P_1^* a & s & gd \\
    0 & 0 & s+\delta-\gamma g d e^{-(\delta_0+s)T}
    \end{array} \right) \\
    &= (s+\mu P_1^* a)(s+\delta_0)(s+\delta-\gamma g d e^{-(\delta_0+s)T})=0.
\end{align*}
    The first two factors give two negative real eigenvalues. Since $N_T<N_{T2}$, it is true that $\gamma g d e^{-\delta_0 T}< \delta$. Set $s=\alpha+i \omega$ for real $\alpha$ and $\omega$ and assume that $\alpha\geq0$. We have that $s+\delta-\gamma g d e^{-(\delta_0+s)T}=0$ if and only if $\alpha + i\omega+\delta=\gamma g d e^{-(\delta_0+\alpha + i\omega)T}$.  Squaring the modulus of both sides we obtain
    \begin{align*}
        (\alpha+\delta)^2 + \omega^2 = (\gamma g d)^2 e^{-2(\delta_0+\alpha)T} < \delta^2,
    \end{align*}
    which implies $\alpha<0$. This contradiction shows that there are no solutions with nonnegative real part. Since all the eigenvalues have negative real parts, $E_1$ is asymptotically stable. \eop
\end{proof}
The global attractivity together with the local stability of $E_1$ for $N_{T1}<N_T<N_{T2}$ shows that $E_1$ is globally asymptotically stable on $D_{N_T}$.

\subsection{Stability of $E_2$: No Delay}

If $m=0$ then the zooplankton are considered mature immediately at birth. In this case there is no delay and system \eqref{model_TDE_2} is just a system of ordinary differential equations. In \cite{kloosterman}, we showed that if $h$ has a negative second derivative (a Type II response) and $f$ has a negative second derivative, which is the case we will consider in the following sections, there is a unique $N_{T3}>N_{T2}$ such that $E_2$ is asymptotically stable if $N_{T2}<N_T<N_{T3}$ and unstable if $N_T>N_{T3}$. This value of $N_{T3}$ is independent of $\delta_0$ and the functional form of $R(P)$, since $m=0$ implies that $\delta_0$ and $R(P)$ do not play a role in system \eqref{model_TDE_2} under the assumption that $(\phi_1,\phi_2,\phi_3)\in D_{N_T}$.

\subsection{Stability of $E_2$: State-Independent Delay}

Consider the case where $R(P)$ is constant. That is, the growth rate of zooplankton is independent of the phytoplankton population. Without loss of generality, we will assume $m$ has been scaled so that $R(P)=1$. In this sense, $m$ has dimensions of time, and represents the age at which zooplankton reaches maturity. In this case, the delay is no longer state-dependent and is fixed at $\tau=m$.

Due to the complicated nature of the resulting characteristic equation, we use numerical methods to study the stability of $E_2$ (both here in the state-independent delay case and in the following section, which deals with the state-dependent delay problem). We use the parameter values in Table \ref{params}.

By setting $s= i\omega$ in the characteristic equation \eqref{characteristic_equation} and allowing the $m$ and $N_T$ parameters to vary, we can find solutions where the eigenvalues have zero real part. These represent critical points where stability might switch. Since the characteristic equation is complex valued, it gives us two equations that need to be satisfied. The equilibrium equations change with changing $m$ and $N_T$, so \eqref{newmodel_eq2} gives us three more equations to solve. We then have six unknowns: $N^*,P^*,Z^*,m,N_T,\omega$. We can then find one-dimensional curves where the five equations are satisfied, using pseudo-arclength continuation \cite{Govaerts}. Note that if $(N^*,P^*,Z^*,m,N_T,\omega)$ solves the five equations then $\omega$ represents the frequency for solutions to system \eqref{model_TDE_2} near the corresponding equilibrium solution $(N^*,P^*,Z^*)$ when the size at which the immature zooplankton reach maturity is $m$ and the total biomass is $N_T$.

\begin{figure}
\begin{center}
\includegraphics[width=0.45 \textwidth]{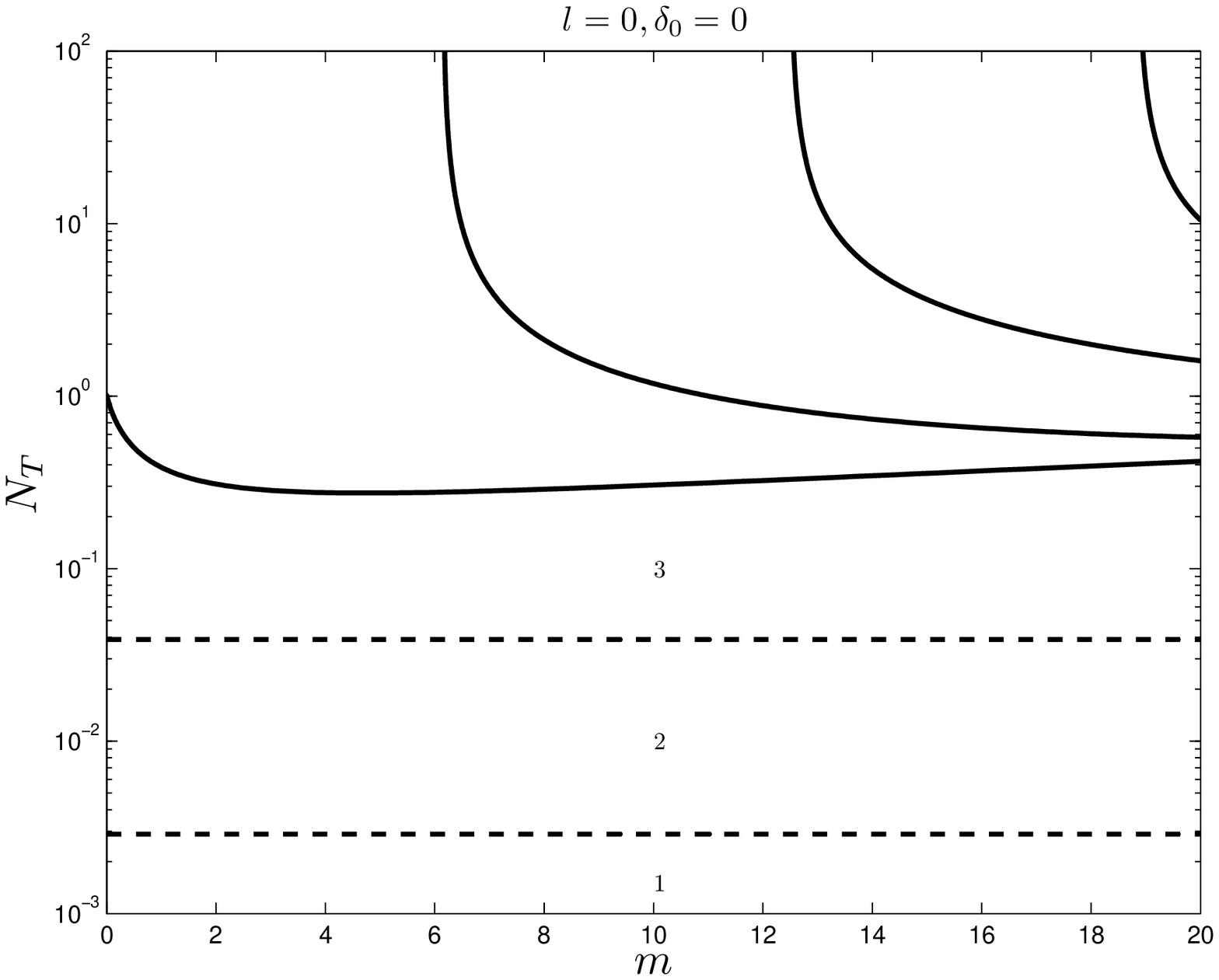}
\includegraphics[width=0.45 \textwidth]{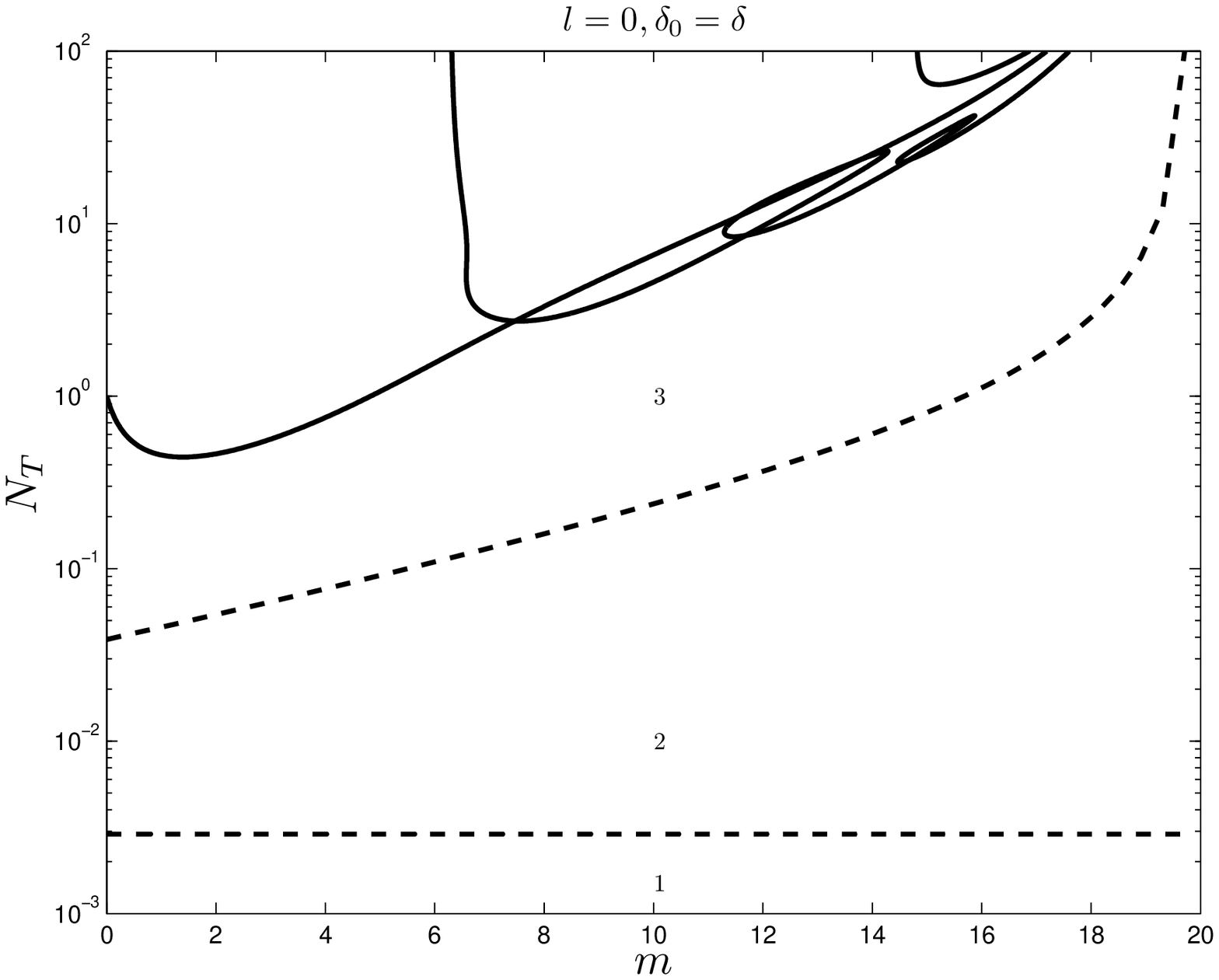}
\end{center}
\caption{The solid curves are where the linearized system has an eigenvalue with zero real part for $R(P)=1$ and $\delta_0=0$ (left) and $\delta_0=\delta$ (right). Region 1 is where $E_1$ and $E_2$ do not exist. Region 2 is where $E_1$ exists, but $E_2$ does not. Region 3 is where $E_1$ and $E_2$ exist. The subset of region 3 under the solid curves is where $E_2$ is asymptotically stable}
\label{const_NT}
\end{figure}

Figure \ref{const_NT} shows regions in the $m-N_T$ plane where the equilibrium solutions exhibit different behaviour.  This was done for two cases: $\delta_0=0$ (left in Figure \ref{const_NT}) and $\delta_0=\delta$ (right in Figure \ref{const_NT}). The former case represents the situation where the immature zooplankton have a zero death rate. That is, we ignore the possibility of them dying before reaching maturity. In region 1 we have that $E_1$ and $E_2$ do not exist. This region is the same, regardless of the value of $\delta_0$. In region 2, $E_1$ exists, but $E_2$ does not. This region does depend on the value of $\delta_0$, as a larger $\delta_0$ requires more biomass for the $E_2$ equilibrium to exist. Region 3 shows where $E_1$ and $E_2$ exist. We note that $R_{\infty}\ln(\gamma g/\delta)/\delta_0 \approx 19.77$, which is an upper limit for $m$ when $\delta_0>0$. The solid curves shows the values of $m$ and $N_T$ where the linearized system has an eigenvalue with zero real part. The subset of region 3 below the solid curves represents values of $m$ and $N_T$ where the $E_2$ equilibrium is asymptotically stable. This follows from the continuity of the eigenvalues with respect to parameter values and the fact that there are no eigenvalues with nonnegative real parts when $m=0$ and $N_{T2}<N_T<N_{T3}$ (between where the upper dotted curve intersects the vertical axis and where the solid curve intersects the vertical axis).

These regions are all fairly independent of $m$ when $\delta_0=0$. That is, if the juvenile zooplankton have a zero death rate, then the required level of maturity has little effect on the stability of the equilibrium solutions. However, when there is positive death rate for the immature zooplankton, the required level of maturity plays a more important role. We see that critical values of total biomass increase with increasing $m$.

\begin{figure}
\begin{center}
\includegraphics[width=0.45 \textwidth]{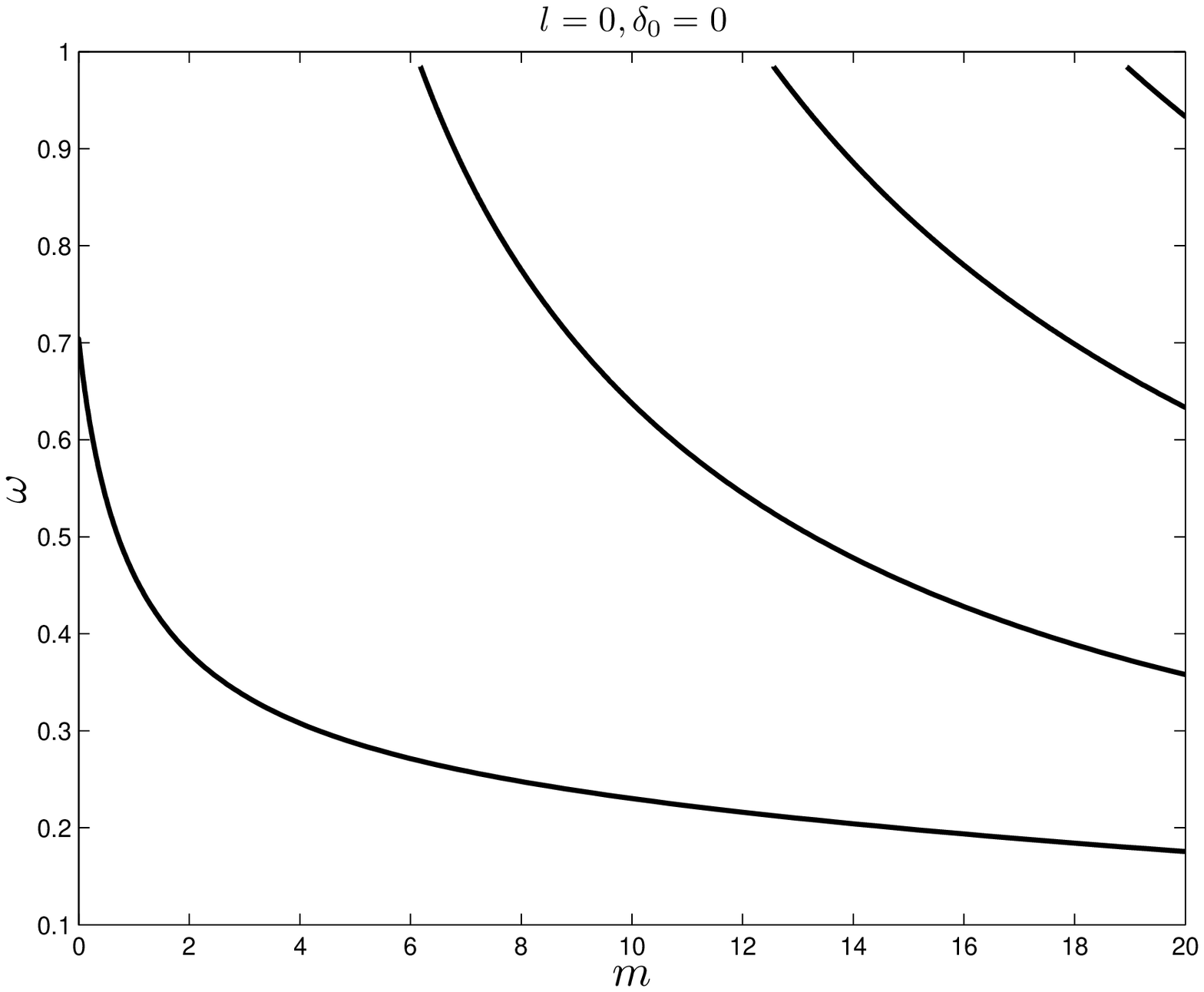}
\includegraphics[width=0.45 \textwidth]{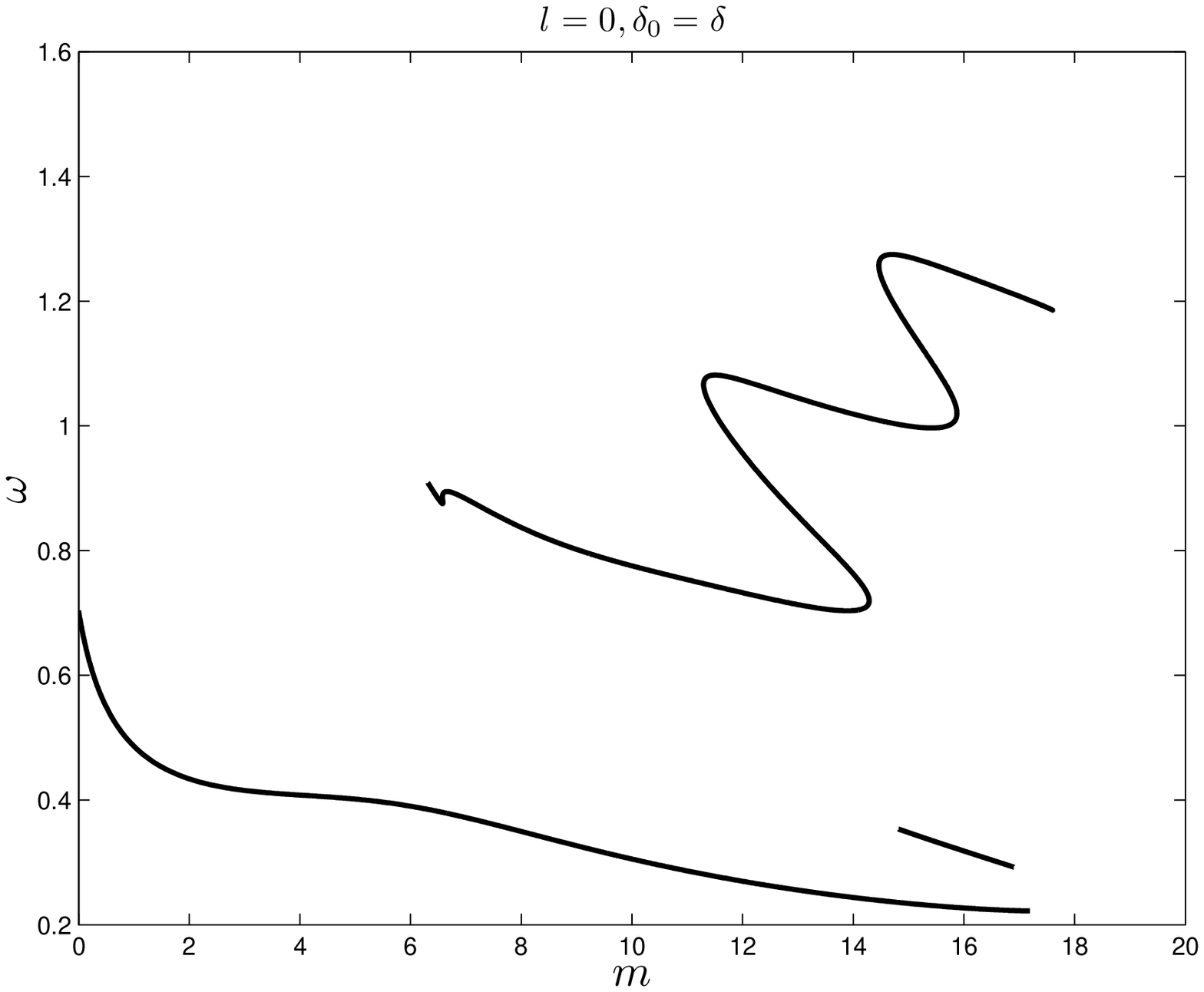}
\end{center}
\caption{The corresponding frequencies to the curves in Figure \ref{const_NT}. The value of $\omega$ is the imaginary part of the eigenvalue with zero real part along these curves}
\label{const_omega}
\end{figure}

Figure \ref{const_omega} shows the frequencies that correspond to the curves in Figure \ref{const_NT}. These are the imaginary parts of the eigenvalues with zero real parts along the curves. For values of $m$ and $N_T$ near the solid curves in Figure \ref{const_NT}, we would expect a slow growth or decay rate in the solution to system \eqref{model_TDE_2} when it is near the equilibrium solution, and for it to have a frequency close to the corresponding value of $\omega$ in Figure \ref{const_omega}.

\subsection{Stability of $E_2$: State-Dependent Delay}

Here we consider the case where $R(P)=\frac{P}{P+l}$ for various values of $l$. Using the same numerical technique as in the previous section, we compute curves in the $m-N_T$ plane where there is an eigenvalue with zero real part. This tells us parameter values where stability can change.

\begin{figure}
\begin{center}
\includegraphics[width=0.45 \textwidth]{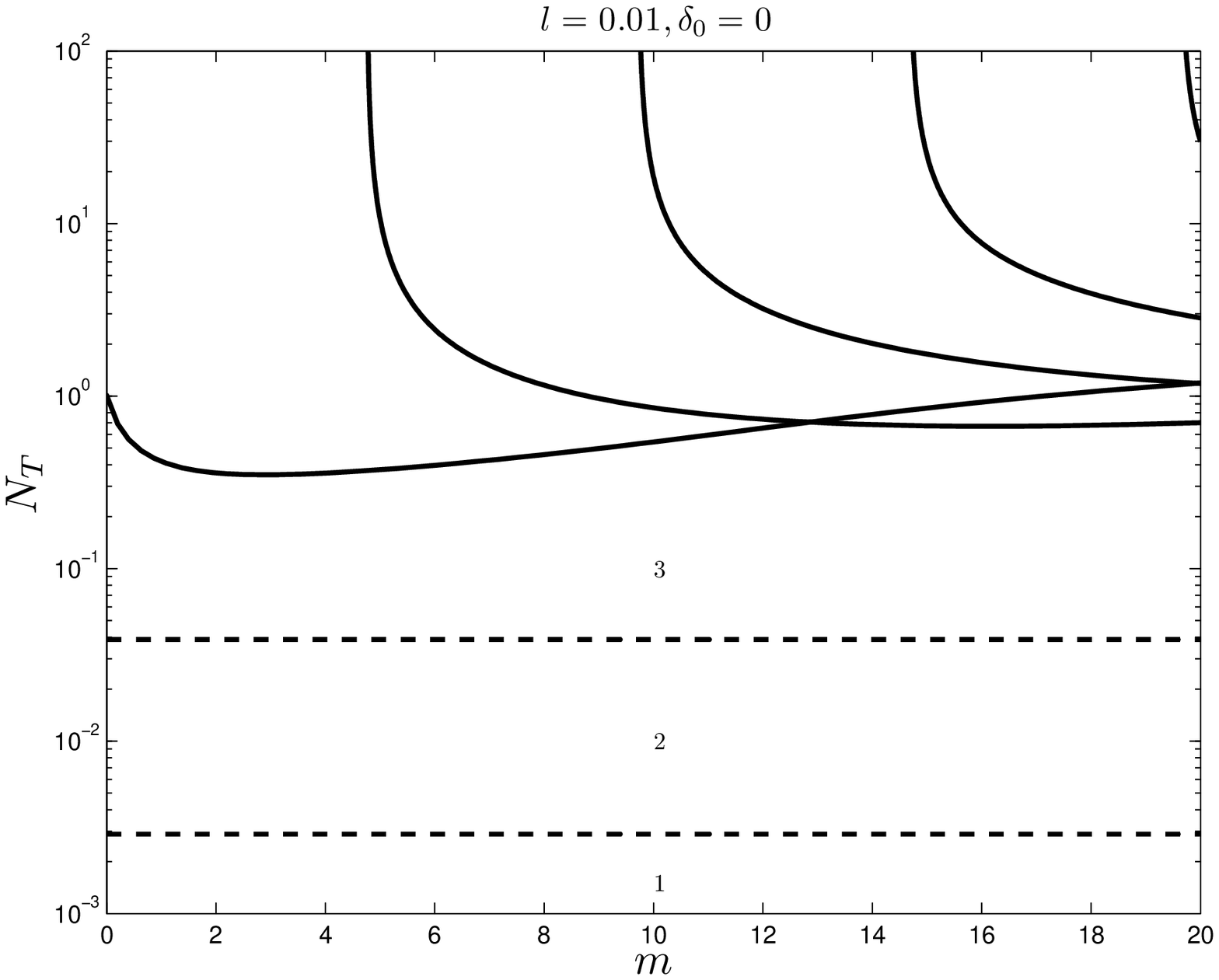}
\includegraphics[width=0.45 \textwidth]{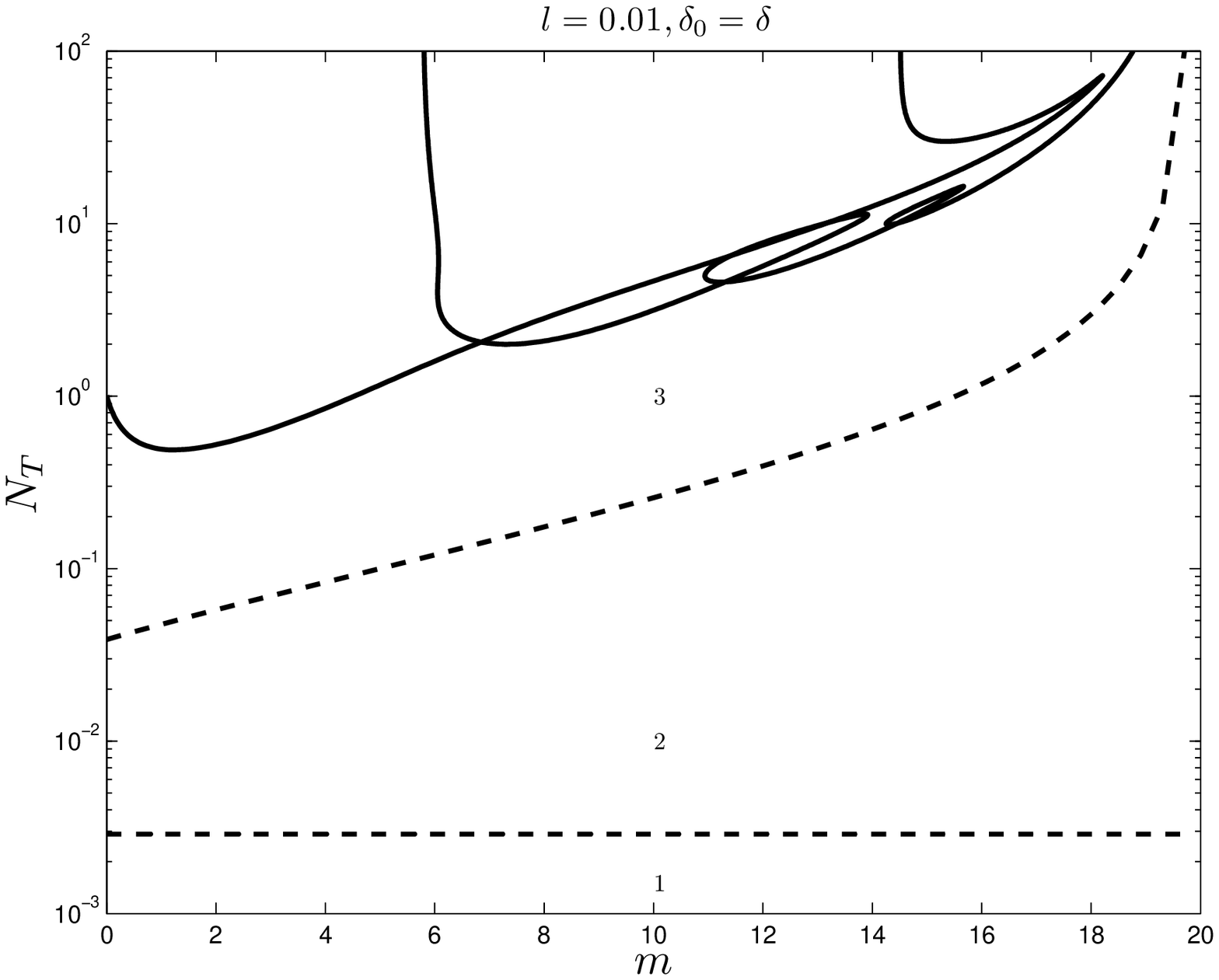}\\
\includegraphics[width=0.45 \textwidth]{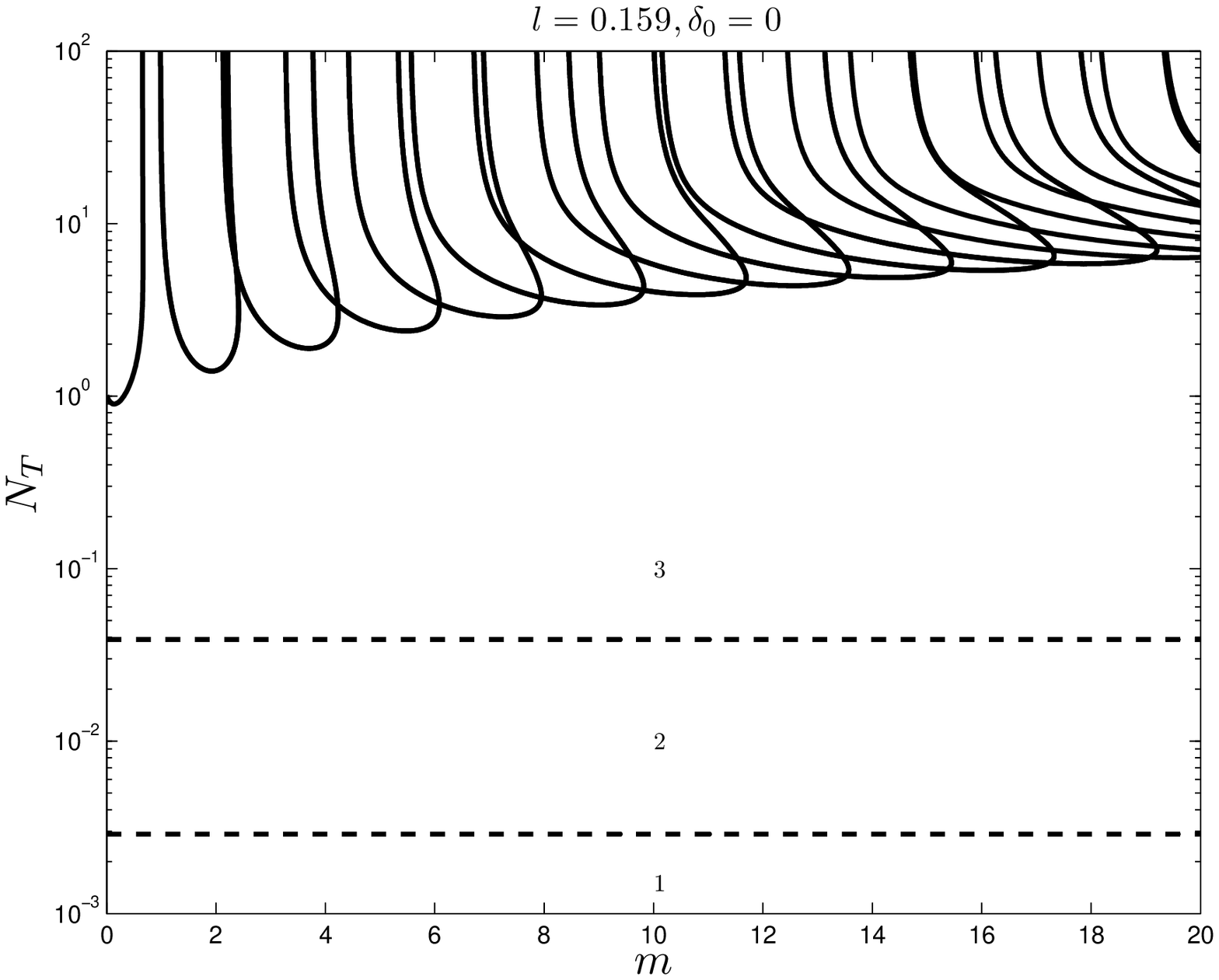}
\includegraphics[width=0.45 \textwidth]{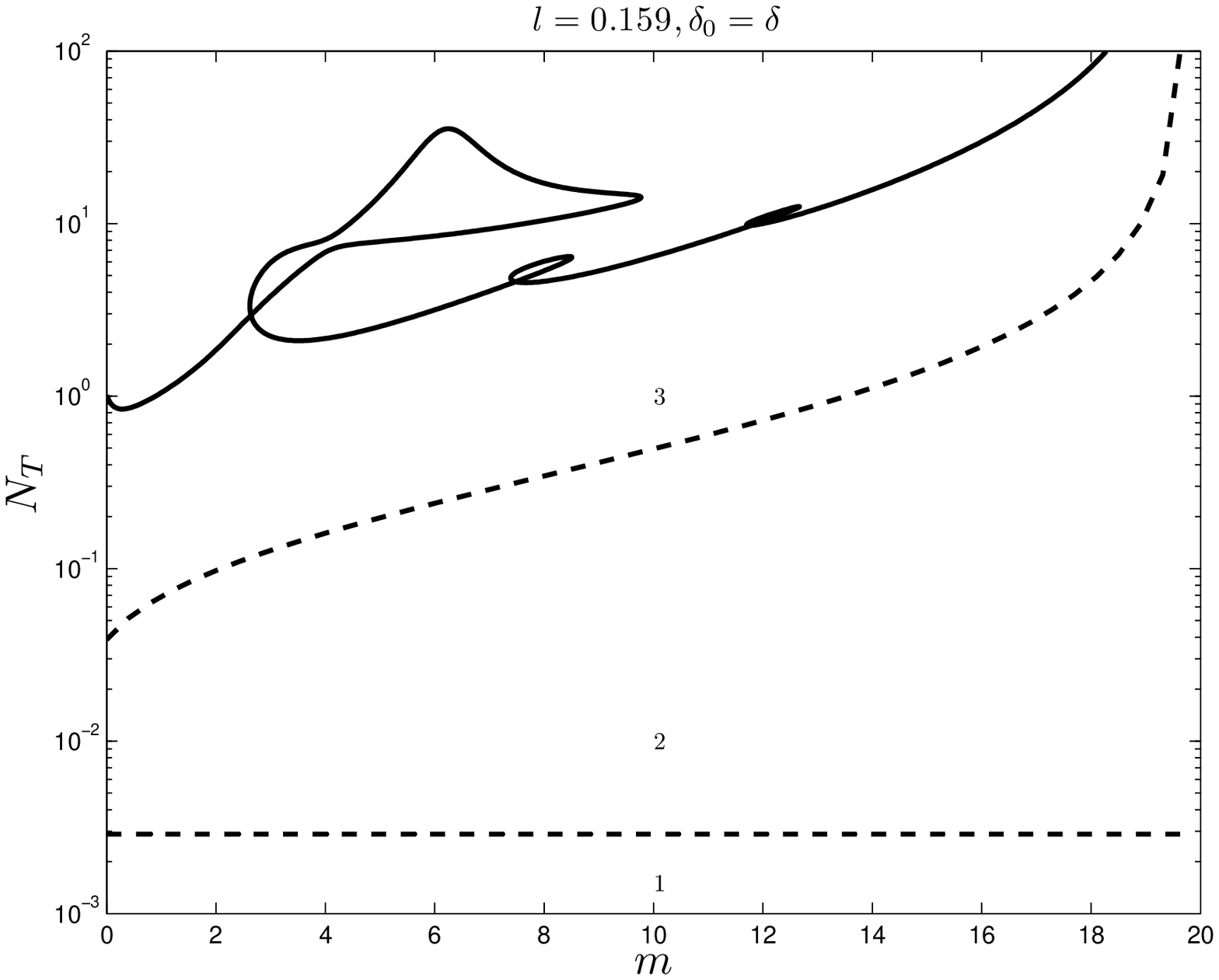}\\
\includegraphics[width=0.45 \textwidth]{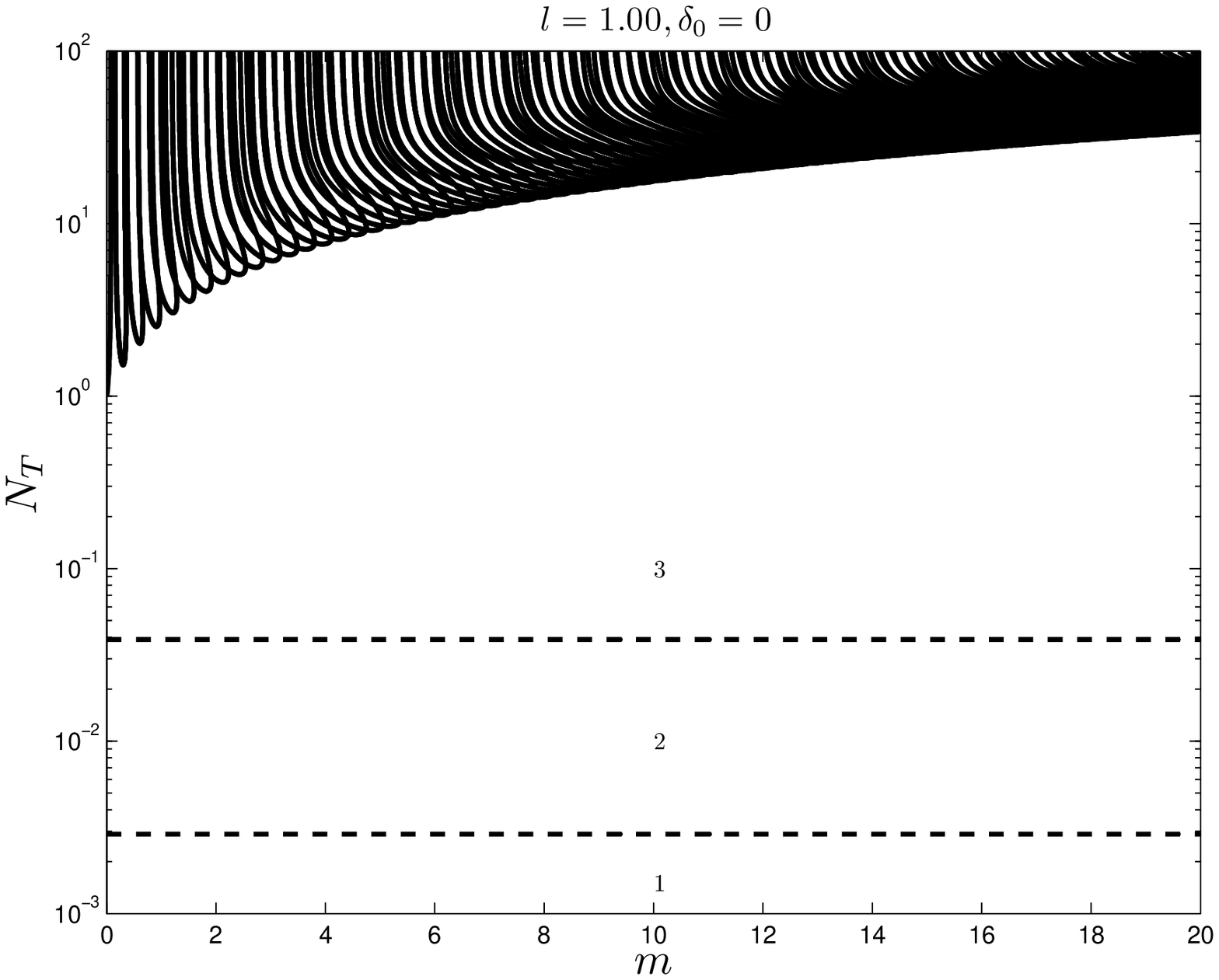}
\includegraphics[width=0.45 \textwidth]{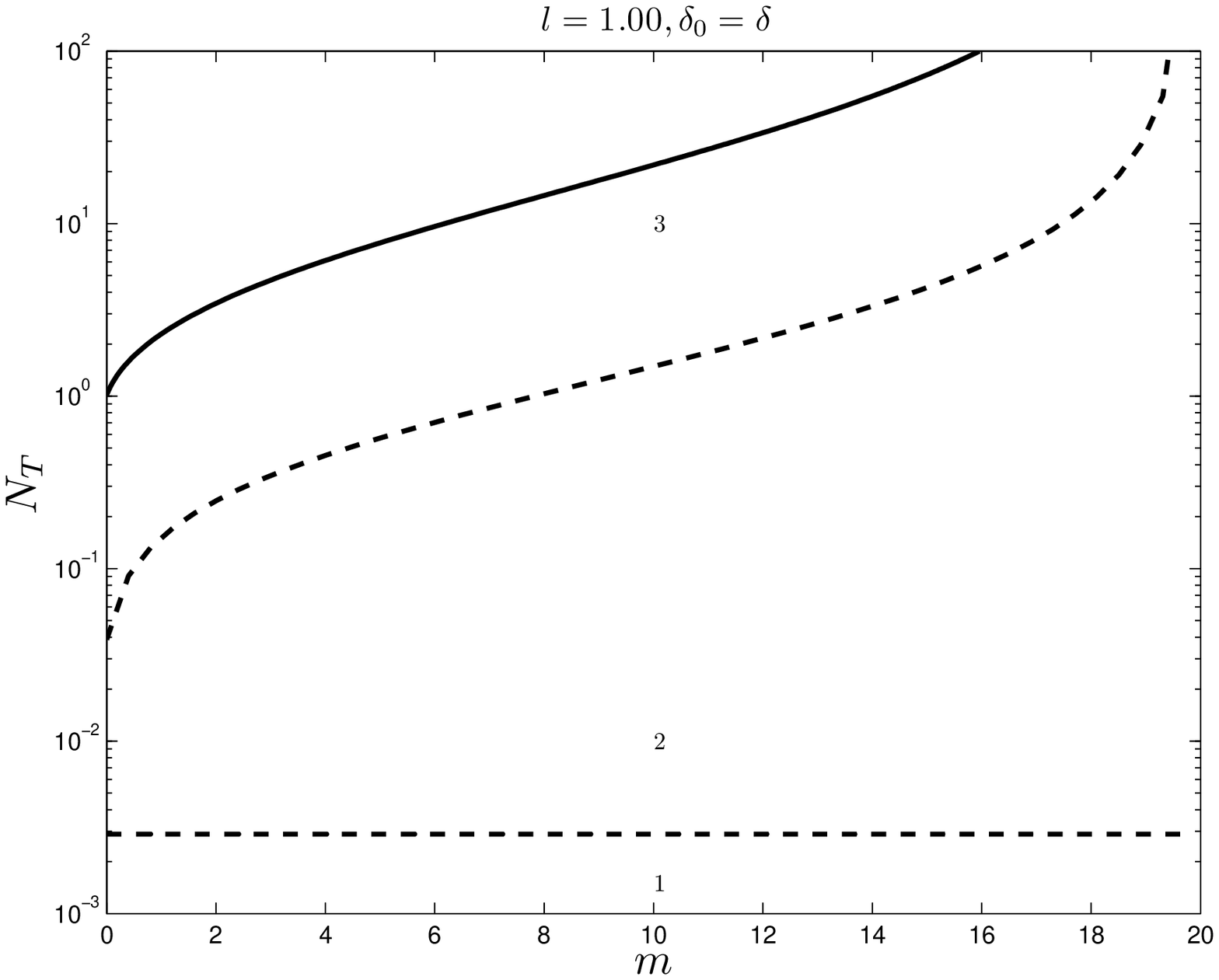}
\end{center}
\caption{The solid curves are where the linearized system has an eigenvalue with zero real part for $R(P)=\frac{P}{P+l}$ for various values of $l$ and $\delta_0=0$ (left) and $\delta_0=\delta$ (right). Region 1 is where $E_1$ and $E_2$ do not exist. Region 2 is where $E_1$ exists, but $E_2$ does not. Region 3 is where $E_1$ and $E_2$ exist. The subset of region 3 under the solid curves is where $E_2$ is asymptotically stable}
\label{sat_NT}
\end{figure}

Figure \ref{sat_NT} shows regions in the $m-N_T$ plane that have different behaviour. This was done for $\delta_0=0$ (left figures), $\delta_0=\delta$ (right figures) and for $l=0.01$ (top figures), $l=0.159$ (figures second from the top), $l=1.00$ (bottom figures). As before, region 1 is where neither $E_1$ nor $E_2$ exist. Region 2 is where $E_1$ exists, but $E_2$ does not. Region 3 is where both $E_1$ and $E_2$ exist. The solid curves are where there is an eigenvalue with zero real part for the system linearized about $E_2$. Hence, the subset of region 3 beneath the solid curves represents the values of $m$ and $N_T$ where the equilibrium solution $E_2$ is asymptotically stable.

We can see that the case where $l=0.01$ is very similar to the case $R(P)=1$, shown in Figure \ref{const_NT}. This is due to the fact that $l$ in this case is small enough relative to $P$ so that $R(P)$ is approximately constant.

When $l$ is increased to $0.159$, its value is no longer small relative to typical values of $P$ (see Figure \ref{equilsolutions}). In this case, the delay has a stronger dependence on the quantity of the phytoplankton.

We see many more curves in the case of $\delta_0=0$. The minimum values of $N_T$ for each of these curves increases slightly as $m$ is increased.  For $l=1.00$ there are about an order of magnitude more curves. The minimum values of $N_T$ for these curves increases more as $m$ is increased, which creates an overall larger region of stability than when $l=0.159$.

For $l=0.159$ and $\delta_0=\delta$ there is a single curve that loops three times. In fact, we observed that $l=0.01$ is a single curve with three loops as well, but we would have to plot it for $N_T$ values much larger than 100 in order to see this. As we vary $l$ from $0.01$ to $0.159$, we see this loop become tighter. As $l$ increases to $1.00$, the loops either vanish or become smaller than what can be detected at this scale. Either way, we get a much simpler region of stability, even though a larger value of $l$ means that the delay has a stronger dependence on the quantity of phytoplankton.

\begin{figure}
\begin{center}
\includegraphics[width=0.45 \textwidth]{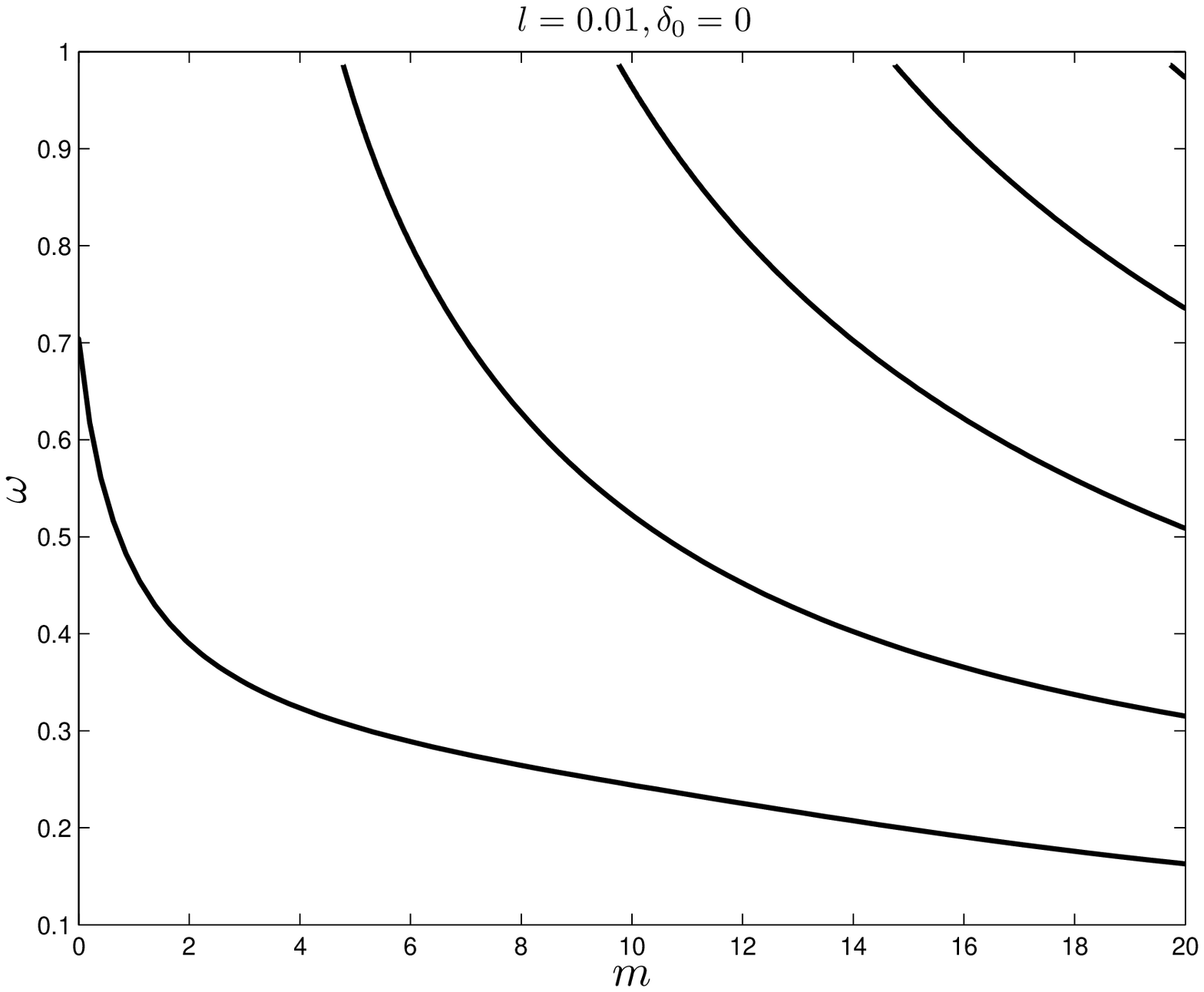}
\includegraphics[width=0.45 \textwidth]{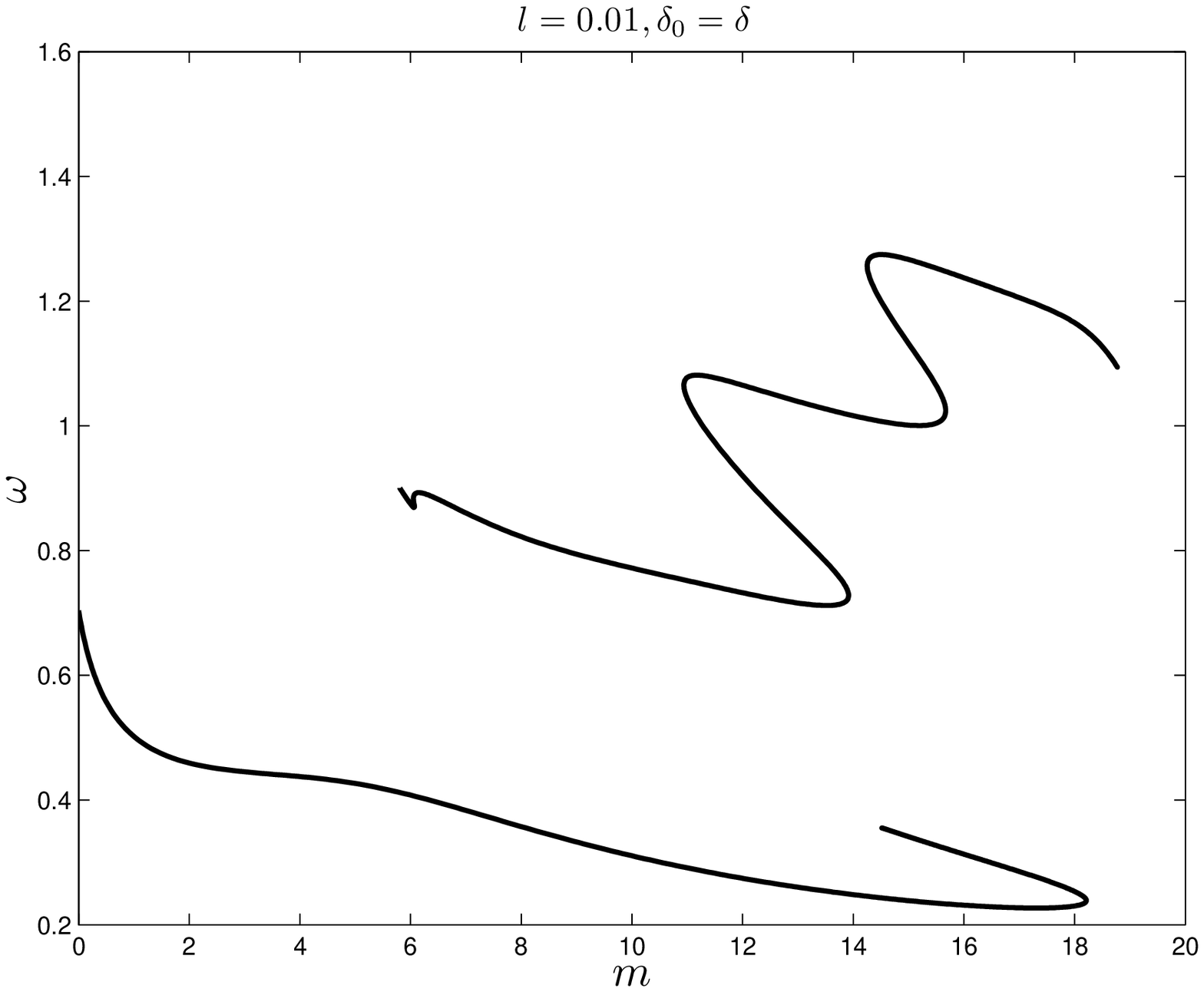}\\
\includegraphics[width=0.45 \textwidth]{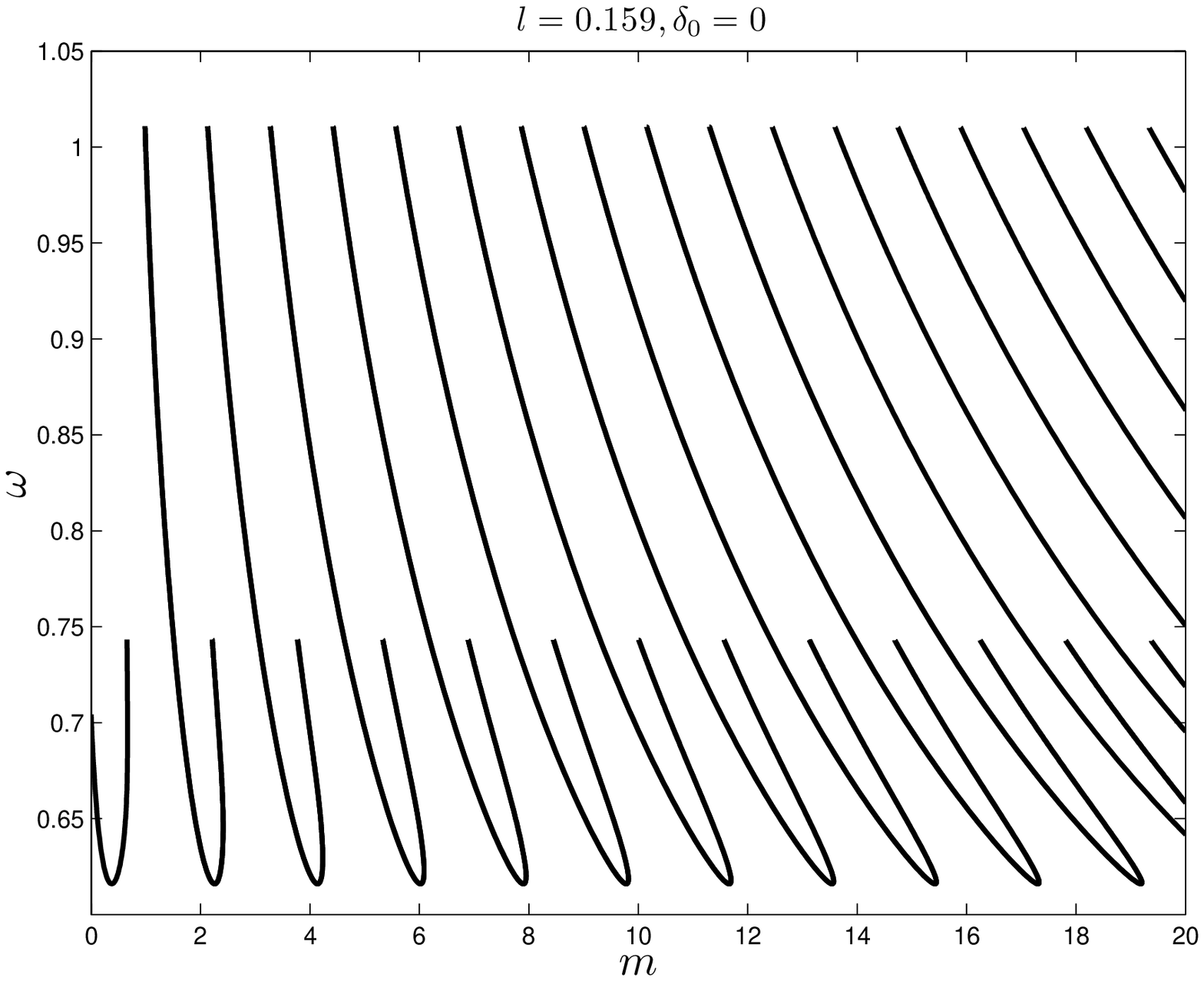}
\includegraphics[width=0.45 \textwidth]{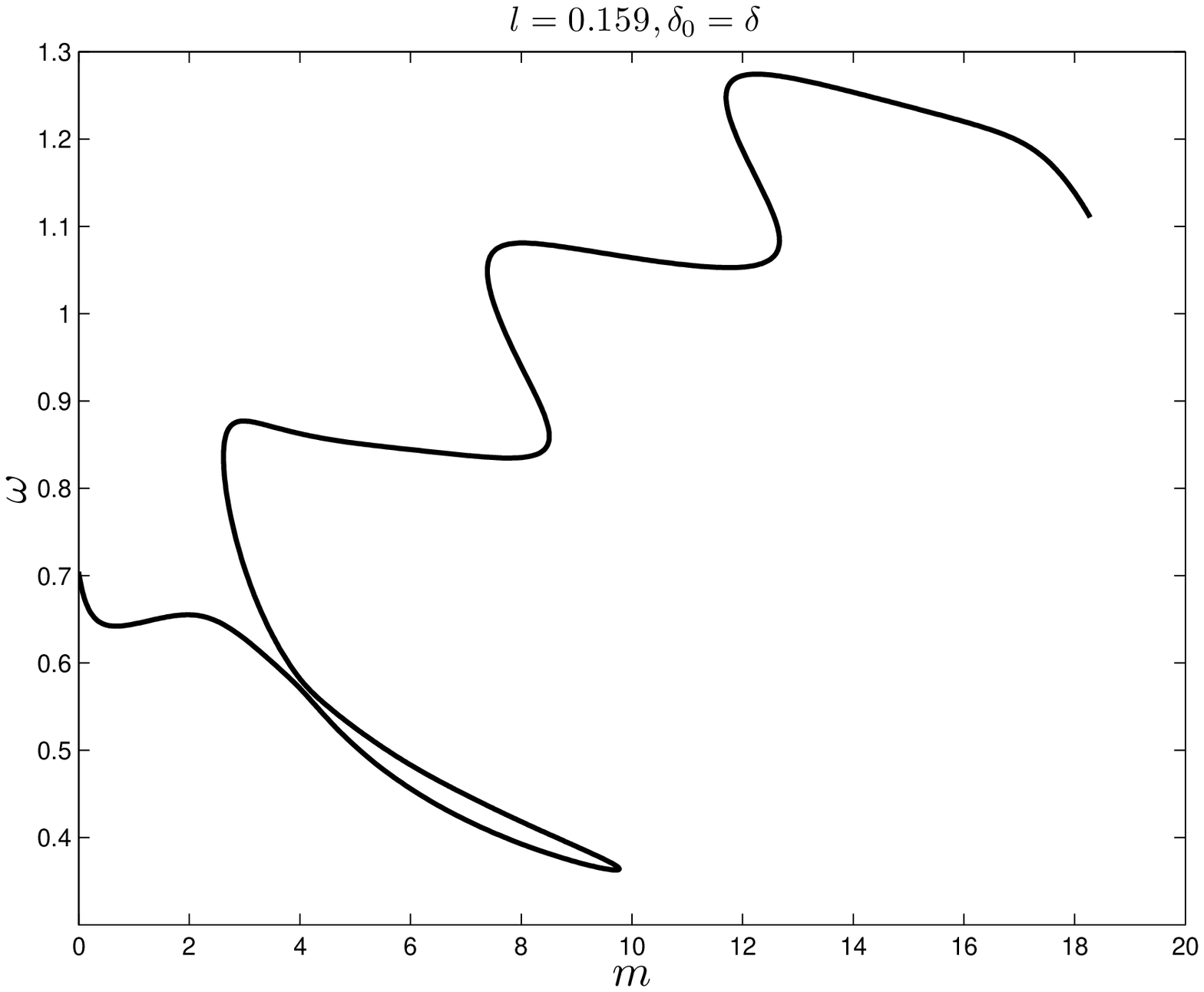}\\
\includegraphics[width=0.45 \textwidth]{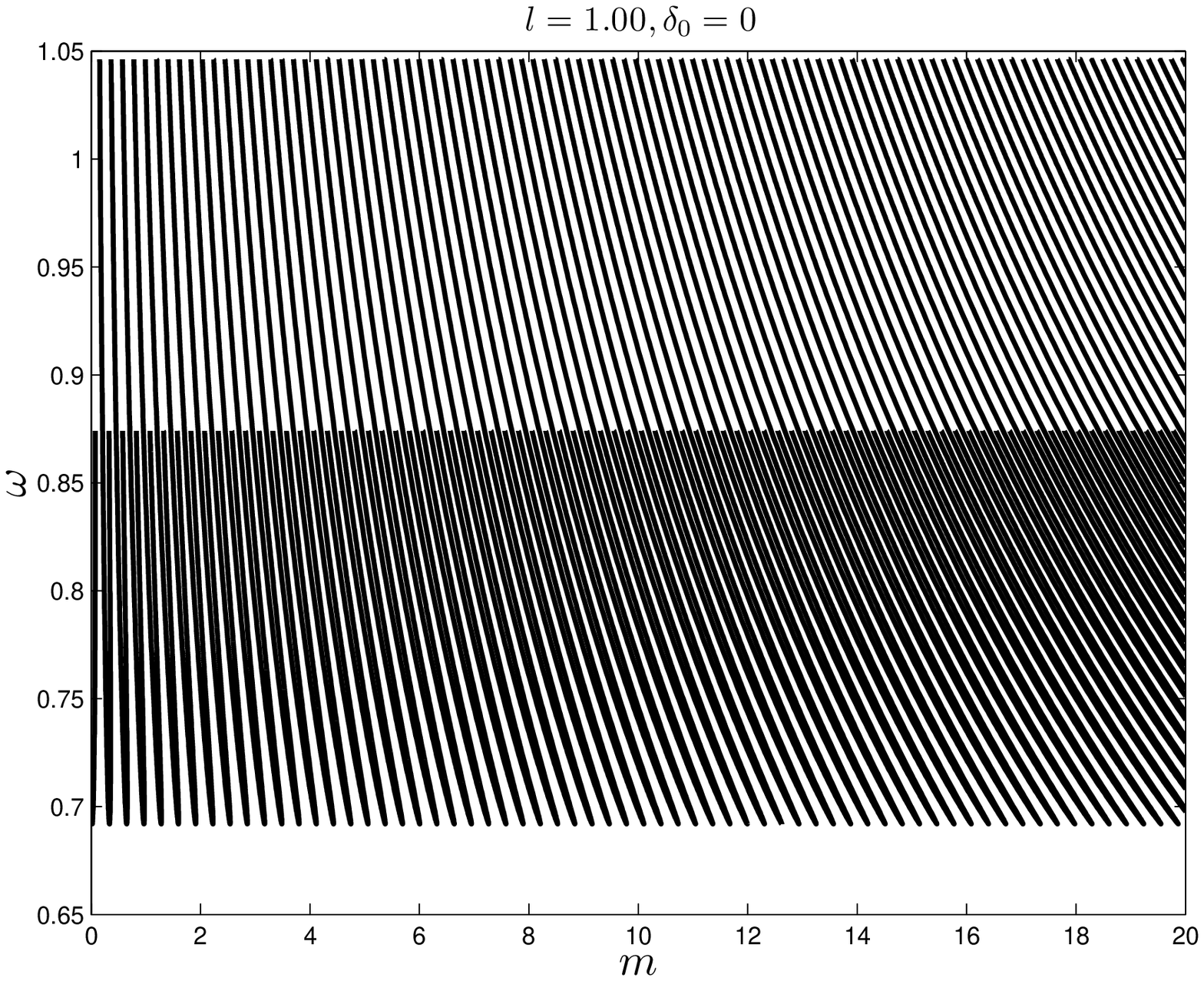}
\includegraphics[width=0.45 \textwidth]{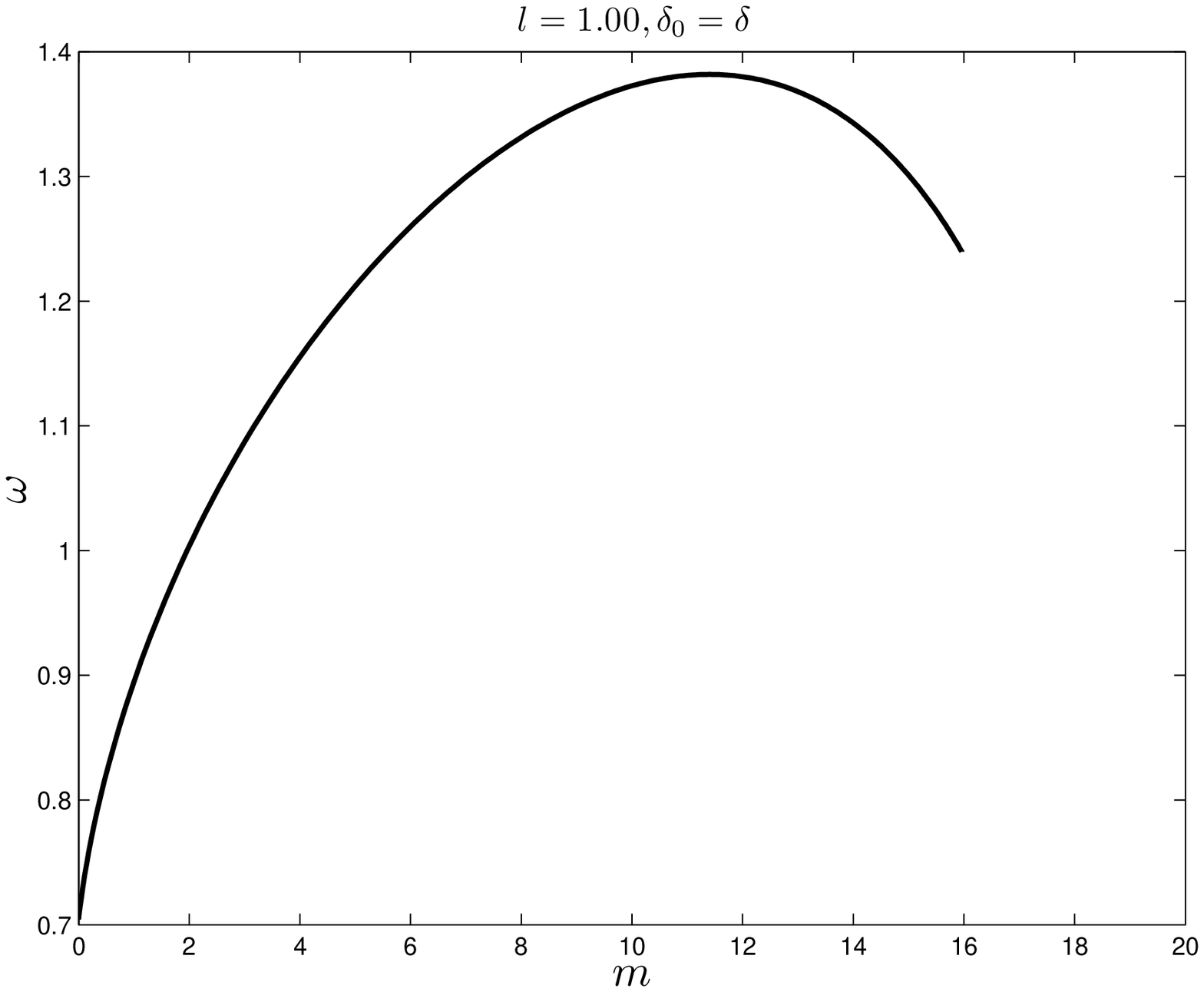}
\end{center}
\caption{The corresponding frequencies to the curves in Figure \ref{const_NT}. The value of $\omega$ is the imaginary part of the eigenvalue with zero real part along these curves}
\label{sat_omega}
\end{figure}

Figure \ref{sat_omega} shows the corresponding frequencies to the curves in Figure \ref{sat_NT}. In this figure, $\omega$ is the imaginary part of the eigenvalues with zero real parts. In the case where $\delta_0=0$, it is true that the characteristic equation is periodic in $m$ in the sense that if the characteristic equation is zero when $m=m_0$, then it is also zero when $m=m_0+2n\pi R(P^*)/\omega$. This is because $m$ only appears in the characteristic equation in the form $e^{i m \omega/R(P^*)}$. This periodicity can be easily seen in Figure \ref{sat_omega} when $\delta_0=0$. Since changing $l$ by an order of magnitude changes $R(P^*)$ by an order of magnitude, we get an order of magnitude more curves for the same range of $m$.

\section{Numerical Simulations}

In order to verify some of the results from the previous section, we performed a series of numerical simulations. Rather than simulating system \eqref{model_TDE_2} directly, we simulated \eqref{model_FDE} and then transformed the resulting time series using equations \eqref{N_transform}-\eqref{Z_transform} and \eqref{time_invert}. System \eqref{model_FDE} was simulated with a second-order method, which included numerically solving the integral \eqref{tauhat} using a trapezoidal rule.

The simulations allowed us to verify the regions of stability for $E_2$ that were computed numerically. For instance, we varied values of $m$ and $N_T$ near the solid curves in Figure \ref{sat_NT}. The initial conditions for the simulation were chosen close to the equilibrium solution $E_2$. For values of $m$ and $N_T$ that were in a stable region, we verified that the time series of the simulation decayed in time. For values of $m$ and $N_T$ that were in an unstable region, we verified that the time series of the simulation grew in time. Many tests were done for various values of $l$ and for $\delta_0=0$ and $\delta_0=\delta$, and no inconsistencies were found. That is, the simulations always agreed with the numerical bifurcation analysis.

\begin{figure}
\begin{center}
\includegraphics[width=0.45 \textwidth]{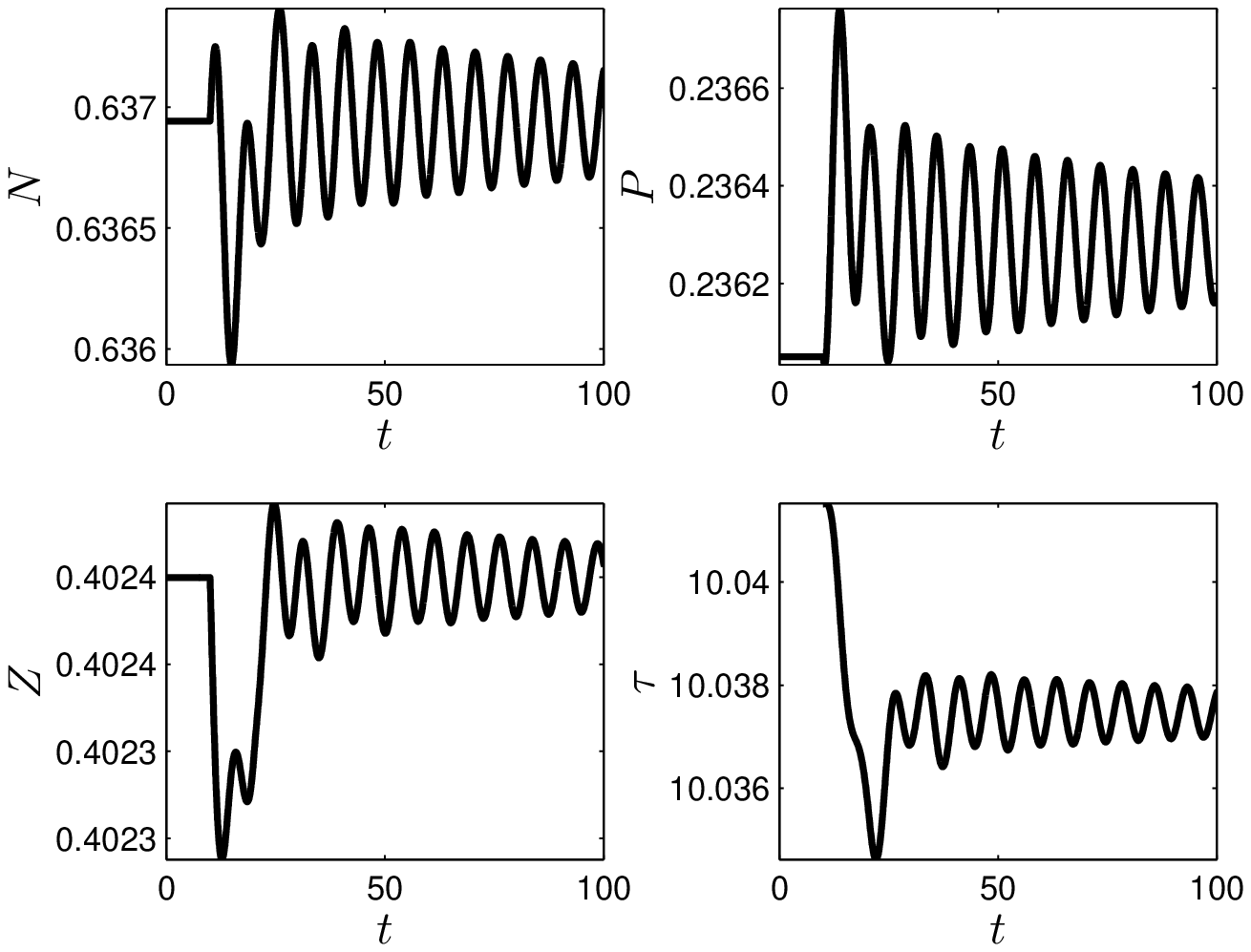}
\includegraphics[width=0.45 \textwidth]{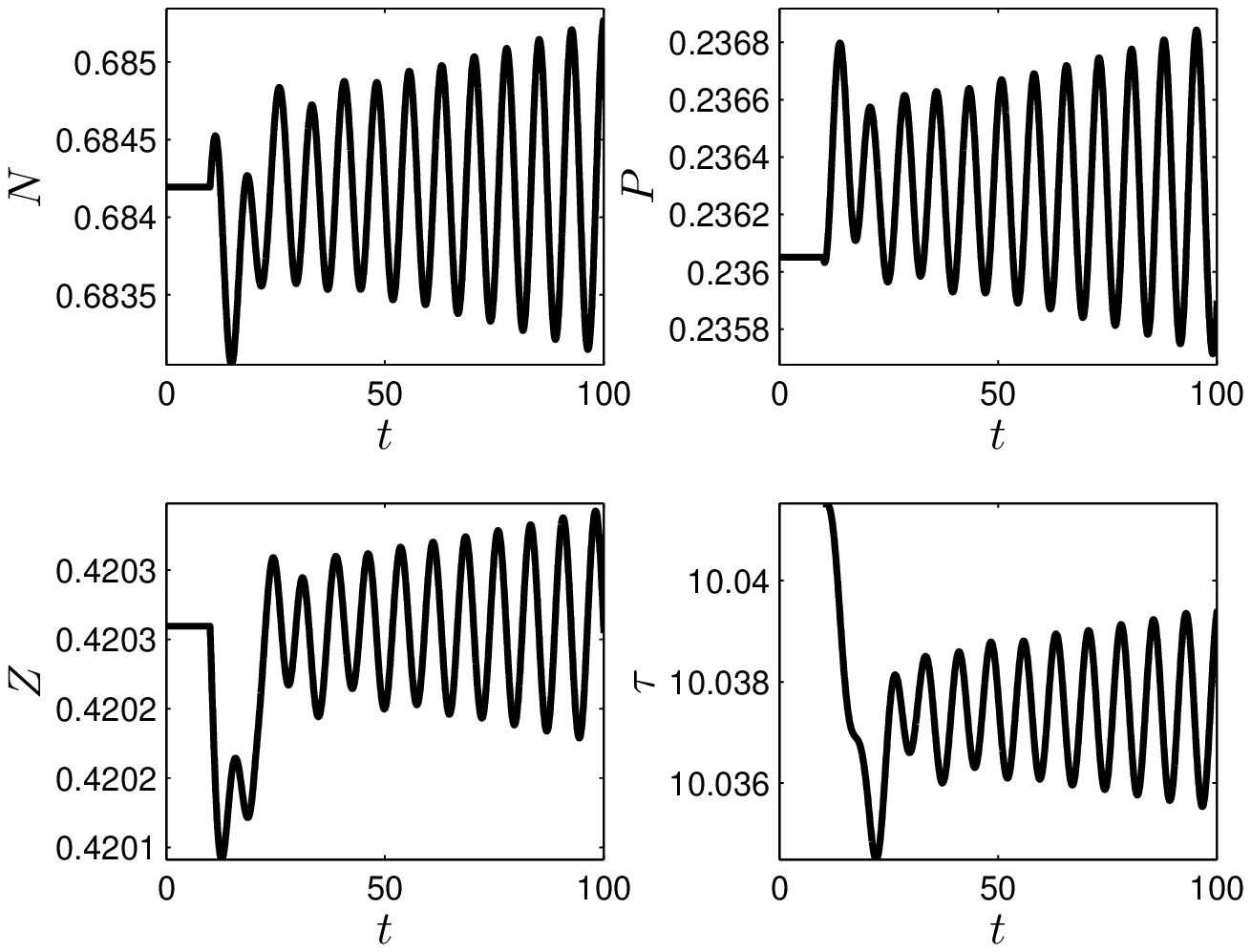}
\end{center}
\caption{The time series of a simulation for $R=\frac{P}{P+0.159}$, $\delta_0=\delta$, $m=6$, and $N_T=10^{0.49}$ (left) and $N_T=10^{0.51}$ (right). The initial conditions were chosen near the equilibrium solution. In the case of $N_T=10^{0.49}$, the simulation suggests that the equilibrium solution is stable, but in the case of $N_T=10^{0.51}$, it suggests that the equilibrium solution is unstable. This agrees with Figure \ref{sat_NT} (on the right, second from the top), which indicated that there is a zero eigenvalue at $m=6$ and $N_T=10^{0.50}$. Also plotted is the time series of the $\tau(m,P_T)$, which is the state-dependent delay}
\label{sim_loc}
\end{figure}

Figure \ref{sim_loc} shows one such verification. For $l=0.159$, $\delta_0=\delta$, and $m=6$, Figure \ref{sat_NT} predicts that $E_2$ should be stable for $N_T=10^{0.49}$ and unstable for $N_T=10^{0.51}$. On the left in Figure \ref{sim_loc}, the simulation for $N_T=10^{0.49}$ indeed suggests that the equilibrium solution is stable. On the right in Figure \ref{sim_loc}, the simulation for $N_T = 10^{0.51}$ suggests that the equilibrium is unstable. While such simulations are by no means a proof of stability, it is still reassuring that they agree with the predictions from the numerical bifurcation analysis.

\begin{figure}
\begin{center}
\includegraphics[width=0.75 \textwidth]{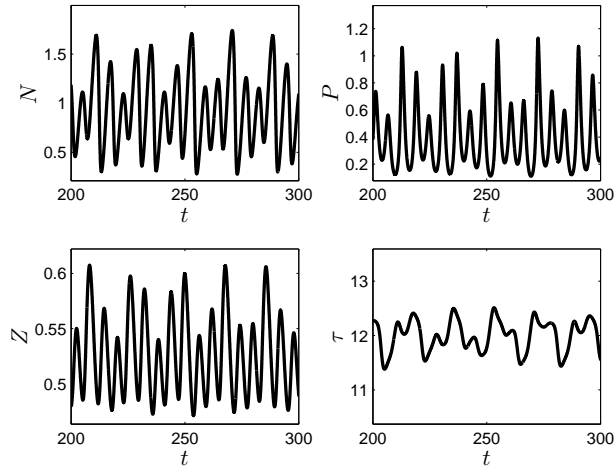}
\end{center}
\caption{A subinterval of the time series of a simulation for $R=\frac{P}{P+0.159}$, $\delta_0=\delta$, $m=8$, and $N_T=10^{0.73}$. The state-dependent delay, $\tau(m,P_t)$ is also shown. This simulation shows the possibility of very irregular and possibly chaotic solutions to system \eqref{model_TDE_2}}
\label{sim_global}
\end{figure}

Figure \ref{sim_global} shows a subinterval of the time series of a simulation for $l=0.159$, $\delta_0=\delta$, $m=8$, and $N_T=10^{0.73}$. For these values Figures \ref{const_NT} and \ref{const_omega} together show that there are about five pairs of eigenvalues with real parts close to zero. Thus, we might expect the solution to exhibit up to five frequencies. In Figure \ref{sim_global} we see that the solution is very irregular, which is not surprising given the predicted spectrum shown in Figure \ref{const_omega}. 

\section{Discussion}

We have looked at a model of a closed planktonic ecosystem that depends on the maturity structure of the immature zooplankton. Using techniques in \cite{smith1991} and \cite{smith1993}, we were able to transform the model into a delay differential equation with a state-dependent threshold-type delay or a delay differential equation with a state-independent delay. Such transformations allowed us to use results readily available from the theory of delay differential equations to study the qualitative features of the model, such as existence and uniqueness of solutions, boundedness, persistence, and stability. For instance, we showed that solutions exist for all time, remain positive, and are bounded, which are desirable features for an ecological model to have.

Being able to represent a system in a variety of ways is very beneficial from a mathematical perspective, as we are able to choose the framework that is most convenient for the situation at hand. For instance, the state-independent delay differential equation was useful for applying well-established theory for functional differential equations with fixed delay, as well as performing numerical simulations. The state-dependent threshold differential equation was convenient for computing the linearization and for arguing qualitative results. The partial differential equation was useful since it represents the model in its most intuitive form, and therefore offers the best framework for interpreting results.

Furthermore, in a practical sense, the different models offer flexibility in what initial data is required for numerical simulations. For instance, in some situations it may be convenient to measure the spectrum of the juvenile zooplankton at a given time, while in other situations it may be more convenient to measure the histories of the phytoplankton and mature zooplankton over a sufficiently long time interval. Since it is possible to transform between the PDE and TDE models, the type of data obtained for the  initial conditions does not necessarily dictate which equations we need to use to simulate the ecosystem.

A key parameter in the study of closed ecosystems is the amount of biomass, which is fixed by the initial conditions. General results include the existence of two critical values of the total biomass. The first is the minimum amount needed to sustain the phytoplankton population, which we have called $N_{T1}$. If the biomass is less than $N_{T1}$, then the phytoplankton and zooplankton both become extinct, as there is not enough biomass to sustain their populations. Conversely, if the biomass is greater than $N_{T1}$, then the phytoplankton do not become extinct. The second critical value of the total biomass, which we have called $N_{T2}$, is the minimum amount needed to sustain the zooplankton population. If the total biomass is less than $N_{T2}$ then the zooplankton population becomes extinct, but does not if the biomass is greater than $N_{T2}$. We have shown that if the total biomass is greater than $N_{T1}$, but less than $N_{T2}$, then the system globally approaches a unique, phytoplankton-only equilibrium solution that depends on the total biomass. Future work might include further study of the system when the biomass is greater than $N_{T2}$, including a formal study of weak/strong (uniform) persistence \cite{butler} and possibly global behaviour of solutions. Here, for $N_T>N_{T2}$, we mainly focused on local stability of the unique equilibrium solution with numerical techniques for a chosen set of parameter values.

In the case when the juvenile zooplankton have a zero mortality rate, ($\delta_0 = 0$), the stability results have a strong dependence on the total biomass ($N_T$) while the required level of maturity ($m$) is less significant. In essence, it seems that the time to maturity does not matter if juvenile zooplankton are not being lost due to mortality. However, when we allow them to have a positive mortality rate ($\delta_0 >0$), the required level of maturity becomes more important with regards to stability. In this case, Figure \ref{sat_NT} suggests that $m$ must increase with $N_T$ in order to maintain stability of the phytoplankton-zooplankton equilibrium.

As well as stable behaviour, the model ecosystem can simple exhibit periodic solutions as well as more complicated dynamics. We saw complicated orbits in many cases where parameters were such that the characteristic equation had multiple pairs of eigenvalues that had real parts close to zero.

Further work may include adding structure to the immature phytoplankton population in a similar way to what was done for the immature zooplankton. In \cite{kloosterman}, we have included a delay in nutrient recycling while ignoring the  structure of the immature zooplankton, but a more complete model would include both of these effects. It may also be worthwhile to see how non-linear closure terms affect the overall behaviour of the system when the zooplankton size-structure is present. An investigation into the role of this closure term in models without such size-structure is given in \cite{edwards2001} and is generally considered to be very significant in determining the dynamics that can occur. Spatial structure can also be added, but similar transformations between partial differential equations and delay differential equations may not be possible in this case.

\bibliography{NPZ_sdd}

\begin{thebibliography}{10}

\bibitem{blythe1984}
S.~P. Blythe, R.~M. Nisbet, and W.~S.~C. Gurney.
\newblock The dynamics of population models with distributed maturation
  periods.
\newblock {\em Theor Popul Biol}, 25(3):289--311, 1984.

\bibitem{bocharov2000}
G.~Bocharov and K.~P. Hadeler.
\newblock Structured population models, conservation laws, and delay equations.
\newblock {\em J Diff Eqn}, 168(1):212--237, 2000.

\bibitem{butler}
G.~Butler, H.~I. Freedman, and P.~Waltman.
\newblock Uniformly persistent systems.
\newblock {\em Proc. Am. Math. Soc.}, 96:425--430, 1986.

\bibitem{cooke1996}
K.~L. Cooke and W.~Huang.
\newblock On the problem of linearization for state-dependent delay
  differential equations.
\newblock {\em Proceedings of the American Mathematical Society}, pages
  1417--1426, 1996.

\bibitem{edwards2001}
A.~M. Edwards.
\newblock Adding detritus to a nutrient-phytoplankton-zooplankton model: a
  dynamical-systems approach.
\newblock {\em J. Plankton Res.}, 23(4):389--413, 2001.

\bibitem{franks2002}
P.~J.~S. Franks.
\newblock {NPZ} models of plankton dynamics: Their construction, coupling to
  physics, and application.
\newblock {\em J. Oceanogr.}, 58:379--387, 2002.

\bibitem{gentleman2008}
W.~C. Gentleman and A.~B. Neuheimer.
\newblock Functional responses and ecosystem dynamics: how clearance rates
  explain the influence of satiation, food-limitation and acclimation.
\newblock {\em J. Plankton Res.}, 30(11):1215--1231, 2008.

\bibitem{gopalsamy}
K.~Gopalsamy.
\newblock {\em Stability and oscillations in delay differential equations of
  population dynamics}.
\newblock Kluwer Academic, Dordrecht, 1992.

\bibitem{Govaerts}
W.~J.~F. Govaerts.
\newblock {\em Numerical Methods for Bifurcations of Dynamical Equilibria}.
\newblock SIAM, Philadelphia, 2000.

\bibitem{gurney1983}
W.~S.~C. Gurney, R.~M. Nisbet, and J.~H. Lawton.
\newblock The systematic formulation of tractable single-species models
  incorporating age structure.
\newblock {\em Journal of Animal Ecology}, 52(2):479--495, 1983.

\bibitem{halebook}
J.~K. Hale and S.~M. Lunel.
\newblock {\em Introduction to Functional Differential Equations}.
\newblock Springer-Verlag, New York, 1993.

\bibitem{holling}
C.~S. Holling.
\newblock The functional response of invertebrate predators to prey density.
\newblock {\em Mem. Entomol. Soc. Can.}, (48):1--86, 1966.

\bibitem{kloosterman}
M.~Kloosterman, S.~A. Campbell, and F.~J. Poulin.
\newblock A closed {N}{P}{Z} model with delayed nutrient recycling.
\newblock {\em J. Math. Biol.}, pages 1--36, 2013.

\bibitem{mccauley}
E.~McCauley, W.~W. Murdoch, and R.~M. Nisbet.
\newblock Growth, reproduction, and mortality of daphnia pulex leydig: life at
  low food.
\newblock {\em Funct. Ecol.}, 4(4):505--514, 1990.

\bibitem{metz}
J.~A.~J. Metz and O.~Diekmann.
\newblock The dynamics of physiologically structured populations.
\newblock {\em Lecture notes in biomathematics}, 68, 1986.

\bibitem{murray}
J.~D. Murray.
\newblock {\em Mathematical biology}.
\newblock Springer-Verlag, Berlin, 1989.

\bibitem{nisbet1983}
R.~M. Nisbet and W.~S.~C. Gurney.
\newblock The systematic formulation of population models for insects with
  dynamically varying instar duration.
\newblock {\em Theoretical Population Biology}, 23(1):114--135, 1983.

\bibitem{poulin}
F.~J. Poulin and P.~J.~S. Franks.
\newblock Size-structured planktonic ecosystems: Constraints, controls and
  assembly instructions.
\newblock {\em J. Plankton Res.}, 32(8):1121--1130, 2010.

\bibitem{smith1991}
H.~L. Smith.
\newblock Threshold delay differential equations are equivalent to standard
  fde’s.
\newblock In {\em Equadiff}, pages 899--904, 1991.

\bibitem{smith1993}
H.~L. Smith.
\newblock Reduction of structured population models to threshold-type delay
  equations and functional differential equations: A case study.
\newblock {\em Mathematical Biosciences}, 113:1--23, 1993.

\bibitem{smith1994}
H.~L. Smith.
\newblock A structured population model and a related functional differential
  equation: Global attractors and uniform persistence.
\newblock {\em Journal of Dynamics and Differential Equations}, 6(1):71--99,
  1994.

\bibitem{smith1995}
H.~L. Smith.
\newblock Equivalent dynamics for a structured population model and a related
  functional differential equation.
\newblock {\em Journal of Mathematics}, 25(1), 1995.

\bibitem{sulsky1989}
D.~Sulsky, R.~R. Vance, and W.~I. Newman.
\newblock Time delays in age-structured populations.
\newblock {\em Journal of theoretical biology}, 141(3):403--422, 1989.

\bibitem{walther}
H.~O. Walther.
\newblock Differential equations with locally bounded delay.
\newblock {\em J. Differential Equations}, 252:3001--3039S, 2012.

\bibitem{waltman1977}
P.~Waltman and E.~Butz.
\newblock A threshold model of antigen-antibody dynamics.
\newblock {\em Journal of Theoretical Biology}, 65(3):499--512, 1977.

\bibitem{wroblewski1988}
J.~S. Wroblewski, J.~L. Sarmiento, and G.~R. Flierl.
\newblock An ocean basin scale model of plankton dynamics in the {N}orth
  {A}tlantic 1. solutions for the climatological oceanographic conditions in
  {M}ay.
\newblock {\em Glob. Biogeochem. Cycles}, 2:199--218, 1988.

\end{thebibliography}
\bibliographystyle{plain}

\end{document}